\documentclass[MSNbibl,number,citesort,seceqn,dvips]{arxbj}
\usepackage{graphicx}

\aid{0}
\volume{19}
\issue{4}
\pubyear{2013}
\firstpage{1268}
\lastpage{1293}
\doi{10.3150/12-BEJSP02} 

\makeatletter
\newcommand{\al}{\alpha}
\newcommand{\be}{\beta}
\newcommand{\cI}{{\mathcal I}}
\newcommand{\cL}{{\mathcal L}}
\newcommand{\cN}{{\mathcal N}}
\newcommand{\cX}{{\mathcal X}}
\newcommand{\De}{\Delta}
\newcommand{\ep}{\varepsilon}
\newcommand{\Ed}{{\mathcal{E}}}
\newcommand{\ga}{\gamma}
\newcommand{\partt}{\partial}
\newcommand{\ra}{\to}
\newcommand{\rai}{\ra\infty}
\newcommand{\si}{\sigma}
\newcommand{\tha}{\theta}
\newcommand{\var}{\operatorname{var}}
\newcommand{\ze}{\zeta}
\newcommand{\ZZ}{\mathbb{Z}}

\newtheorem{them}{Theorem}

\setattribute{abstract}    {skip} {40\p@}
\setattribute{keyword}     {skip} {16\p@}
\setattribute{frontmatter} {skip} {0\p@ plus 3\p@ minus 3\p@}

\makeatother

\begin{document}
\begin{frontmatter}

\title{Normal approximation and smoothness for sums of means of lattice-valued random variables}
\runtitle{Sums of independent lattice-valued random
variables}\vspace*{12pt}

\begin{aug}
\author[1]{\fnms{Geoffrey} \snm{Decrouez}\thanksref{1}\ead[label=e1]{dgg@unimelb.edu.au}} \and
\author[2]{\fnms{Peter} \snm{Hall}\corref{}\thanksref{2}\ead[label=e2]{halpstat@ms.unimelb.edu.au}}
\runauthor{G. Decrouez and P. Hall} 
\address[1]{Department of Mathematics and Statistics, The University
of Melbourne, VIC 3010, Australia.\\
\printead{e1}}
\address[2]{Department of Mathematics and Statistics, The University
of Melbourne, VIC 3010, Australia and Department of Statistics,
University of California Davis, Davis, CA 95616, USA.\\
\printead{e2}}
\end{aug}


%
\begin{abstract}
Motivated by a problem arising when analysing data from quarantine
searches, we explore properties of distributions of sums of independent
means of independent lattice-valued random variables. The aim is to
determine the extent to which approximations to those sums require
continuity corrections. We show that, in cases where there are only two
different means, the main effects of distribution smoothness can be
understood in terms of the ratio $\rho_{12}=(e_2 n_1)/(e_1 n_2)$,
where $e_1$ and $e_2$ are the respective maximal lattice edge widths of
the two populations, and $n_1$ and $n_2$ are the respective sample
sizes used to compute the means. If $\rho_{12}$ converges to an
irrational number, or converges sufficiently slowly to a rational
number; and in a number of other cases too, for example those where
$\rho_{12}$ does not converge; the effects of the discontinuity of
lattice distributions are of smaller order than the effects of
skewness. However, in other instances, for example where $\rho_{12}$
converges relatively quickly to a rational number, the effects of
discontinuity and skewness are of the same size. We also treat
higher-order properties, arguing that cases where $\rho_{12}$
converges to an algebraic irrational number can be less prone to suffer
the effects of discontinuity than cases where the limiting irrational
is transcendental. These results are extended to the case of three or
more different means, and also to problems where distributions are
estimated using the bootstrap. The results have practical
interpretation in terms of the accuracy of inference for, among other
quantities, the sum or difference of binomial proportions.\looseness=1
\end{abstract}

%
\begin{keyword}
\kwd{algebraic irrational number}
\kwd{bootstrap}
\kwd{confidence interval}
\kwd{continuity correction}
\kwd{difference of binomial proportions}
\kwd{discontinuity}
\kwd{irrational number}
\kwd{sum of binomial proportions}
\kwd{transcendental number}
\end{keyword}

\end{frontmatter}

\section{Introduction}\label{sec1}

\subsection{Background: The case of a single sample mean}\label{sec1.1}

Let ${\hat\tha}$ denote a statistical estimator of an unknown
quantity $\tha
$, and assume that ${\hat\tha}-\tha$ is asymptotically normally distributed
with zero mean and variance $n^{-1}\si^2$, where $n$ is a measure of
the sample size from which ${\hat\tha}$ was computed. In particular, the
statistic $T=n^{1/2}({\hat\tha}-\tha)/\si$ is asymptotically normal
$\mathrm{N}(0,1)$. Under additional assumptions an Edgeworth expansion of the
distribution of $T$ can generally be formulated, having the form
%
%
\begin{equation}
\label{eq1.1} P(T\leq x)=\Phi(x)+n^{-1/2}P(x) \phi(x)+o
\bigl(n^{-1/2} \bigr),
\end{equation}
uniformly in $x$, where $\Phi$ and $\phi$ are the standard normal
distribution and density functions, respectively, and $P$ is an even,
quadratic polynomial.

For example, if ${\hat\tha}$ denotes the mean of a sample of size $n$
from a
population with mean~$\tha$, and if data from the population have
finite third moment and a nonlattice distribution, then (\ref{eq1.1}) holds
with $P(x)={\frac{1}{6}}\be(1-x^2)$, where $\be$ is the standardised skewness
of the population. On the other hand, still in the case of the mean of
a population with finite third moment, if the population is lattice
then an extra, discontinuous term has to be added to (\ref{eq1.1}).

This extra term reflects the discrete continuity correction that
statisticians are often obliged to introduce when approximating a
lattice distribution, for example the binomial distribution or the
distribution of a sum of Poisson variates, using the smooth normal distribution:
%
%
\begin{equation}
\label{eq1.2} P(T\leq x)=\Phi(x)+n^{-1/2}P(x) \phi(x)
+n^{-1/2}d_n(x) \phi(x)+o \bigl(n^{-1/2} \bigr),
\end{equation}
where $d_n(x)=(e_0/\si) \psi_n(x)$ denotes the discontinuous term in
the Edgeworth expansion, $e_0$ is the maximal span of the lattice, $\si^2$ is the population variance,
\[
\psi_n(x)=\psi \bigl\{(x-\xi_n) \si n^{1/2}
/e_0 \bigr\}, \qquad \xi_n= \bigl(e_0 /\si
n^{1/2} \bigr) \bigl\{{\tfrac{1}{2}}-\psi(n x_0/e_0)
\bigr\},
\]
$\psi(x)=\lfloor x\rfloor-x+{\frac{1}{2}}$, $\lfloor x\rfloor$ is the
largest integer
not strictly exceeding $x$, and it is assumed that all points of
support in the distribution (of which $\tha$ is the mean) have the
form $x_0+\nu e_0$ for an integer $\nu$.

\subsection{Contributions of this paper}\label{sec1.2}

We show that in multi-sample problems, where ${\hat\tha}$ is a sum of
several independent means, the discontinuous term can be ignored if
sample sizes are chosen judiciously. For example, if there are just two
sample sizes (as in the case of a sum or difference of two binomial
proportions), and if the lattice edge widths are identical (this
simplifies our discussion here, but is not essential), then it is
sufficient that the ratio of the two sample sizes converge to an
irrational number, or converge sufficiently slowly to a rational
number. These results are corollaries of Theorem~\ref{the1} in
Section~\ref{sec2.1}, and
they and other properties are discussed in
Section~\ref{sec2.2}.\vadjust{\eject}

We also show that the discontinuous term can be replaced by
$O(n^{\delta-1})$, for all $\delta>0$, provided that four moments are
finite and the ratio of the two sample sizes converges sufficiently
quickly to an irrational number of ``type'' 1. (See Section~\ref
{sec2.3} for a
definition of the type of an irrational number.) More generally, we
explore the effect that type has on the size of the discontinuous term.
Theorem~\ref{the2} also gives an explicit formula for the
discontinuous term, up
to a remainder of order $n^{-1}$. Sections~\ref{sec2.5} and~\ref
{sec2.6} show how to bound
the discontinuous remainder term, for two different approaches to
defining that quantity, and show how the effects of irrationals of
different types can be teased from the remainder. Applications to the
bootstrap are straightforward, and in fact Section~\ref{sec2.7}
outlines a
bootstrap version of Theorem~\ref{the2} and discusses its implications.

We do not treat in any detail cases where the differences between two
lattice distributions arise mainly in terms of their centres, rather
than their lattice edge widths. For example, if two independent sample
means ${\bar X}_j$, for $j=1,2$, are respectively averages of $n_j$
independent variables and are distributed on lattices $x_j+\nu n_j^{-1}
e_j$; and if the difference $x_1-x_2$ between the lattice translations
equals an irrational multiple of the ratio $\rho_{12}=(e_2 n_1)/(e_1
n_2)$; then the distribution of ${\bar X}_1+{\bar X}_2$ is non-lattice. While
this problem and its implications are of mathematical interest, they do
not enjoy the practical motivation of problems where, say, $x_1=x_2$
and $\rho_{12}$ can be almost arbitrary. For example, $x_1=x_2$ in the
problem of constructing confidence intervals for the sum or difference
of two binomial probabilities, based on samples of unequal size.
Therefore, we address cases where the focus of attention is $\rho_{12}$ rather than $x_1-x_2$. Differences between lattice centres are
permitted by our regularity conditions, but their role is not treated
in detail.

\subsection{Practical motivation}\label{sec1.3}

The extra term in (\ref{eq1.2}), relative to (\ref{eq1.1}), is of
significant interest
to a practitioner, since it causes significant inaccuracy when the
central limit theorem is used to approximate the distribution of $T$.
The presence of this extra term motivates the continuity correction,
and also the fiducial approach taken by Clopper and Pearson \cite
{Clo34} and Sterne \cite{Ste54} to estimating a binomial proportion,
as well as a large, more recent literature discussing methodology for
solving problems such as constructing confidence intervals for the
difference or sum of two binomial proportions. See, for example, Hall
\cite{Hal82}, Duffy and Santer \cite{Duf87}, Lee et~al.~\cite{Lee97},
Agresti and Caffo \cite{Agr00}, Brown et al. \cite{Bro01,Bro02}, Zhou
et al. \cite{Zho01}, Price and Bonnett \cite{Pri04}, Brown and Li
\cite{Bro05}, Borkowf \cite{Bor06}, Roths and Tebbs \cite{Rot06},
Wang \cite{Wan10} and Zieli\'nski \cite{Zie10}.

The practical motivation for the work described in this paper came from
data acquired during quarantine searches, where the construction of
confidence intervals for the sum, rather than difference, of two
binomial proportions was of interest. In detail, shipping containers
arriving at a frontier contained a certain number, $N$ say, of
consignments. Some of the consignments might be clean, but others could
contain pests which needed to be detected and removed to prevent their
introduction to the environment. To reduce the costs associated with
inspection, quarantine services usually inspect only $n_1<N$
consignments. Consignments are assumed to be contaminated with
probability $p_1$, and the number, $n_1 {\bar X}_1$ say, of contaminated
consignments found after routine (but incomplete) inspection of the
items in each of $n_1$ consignments is assumed to follow a binomial
distribution. Contaminated consignments are then ``cleaned,'' and the
members of a subsample of $n_2$ of them are reinspected. (The data
gathered in this way comprise a ``leakage survey.'') The number of
items, $n_2 {\bar X}_2$, still found contaminated (for example,
contaminated with a different kind of pest) are assumed to follow a
binomial distribution with parameters $n_2$ and $p_2$, and typically it
is argued that ${\bar X}_1$ and ${\bar X}_2$ are independent. An
estimator of
the proportion of items that pass through this inspection process, and
are still contaminated, is given by
%
%
\begin{equation}
\label{eq1.3} {\bar X}_1 \biggl(1-\frac{n_1}{ N} \biggr)+{\bar
X}_2 \biggl(1-\frac
{n_2}{ N} \biggr) ,
\end{equation}
which can be viewed as a sum of means of lattice-valued random
variables where the lattice edge lengths are $e_j=1-N^{-1}n_j$ for $j=1,2$.

The quarantine inspection service aims to develop a strategy for
choosing consignments, and items, to inspect. This reduces the
associated costs, and minimises, to at least some extent, the number of
contaminated items that cross the border. The performance of such a
strategy is assessed, by the quarantine service, using a variety of
statistics based on sums of binomials; (\ref{eq1.3}) is just one example.
Quarantine services are usually interested in providing confidence
intervals as well as point estimators, and hence there is significant
interest in estimating the distributions of statistics such as that at~(\ref{eq1.3}).

\section{Main results}\label{sec2}

\subsection{Edgeworth expansions with remainder equal to $o(n^{-1/2})$}\label{sec2.1}

Let $X_{ji}$, for $1\leq i\leq n_j$ and $j=1,\ldots,k$, denote
independent random variables. Assume that each $X_{ji}$ has a
nondegenerate lattice distribution, depending on $j$ but not on $i$ and
with maximal lattice edge width $e_j$ and finite third moment. Suppose
too that $k\geq2$. Put ${\bar X}_j=n_j^{-1}\sum_iX_{ji}$, $\mu_j=E(X_{ji})$,
$\si_j^2=\var(X_{ji})$ and
%
\begin{equation}
\label{eq2.1} S=\sum_{j=1}^k {\bar
X}_j.
\end{equation}
The model (\ref{eq2.1}) includes cases of apparently greater
generality, for
example where signed weights are incorporated in the series, since the
absolute values of the weights can be incorporated into (\ref{eq2.1}) by
modifying the lattice edge widths, and negative signs can be addressed
by reflecting the summand distributions.

Since third moments are finite then, if the distributions of
$X_{11},\ldots,X_{k1}$ were to satisfy a smoothness condition, such as
that of Cram\'er, we could express the distribution of $S$ in a
one-term Edgeworth expansion:
%
%
\begin{equation}
\label{eq2.2} P \biggl\{\frac{S-E(S)}{(\var S)^{1/2}}\leq x \biggr\} =\Phi(x)+n^{-1/2}{
\frac{1}{6}}\be \bigl(1-x^2 \bigr) \phi(x)+o
\bigl(n^{-1/2} \bigr) ,
\end{equation}
where we take $n=n_1+\cdots+n_k$ to be the asymptotic parameter, and
%
%
\begin{equation}
\label{eq2.3} \be=\be(n)=\frac{n^{1/2}E(S-ES)^3}{(\var S)^{3/2}} =\frac{n^{1/2}\sum_jn_j^{-2}E(X_{j1}-EX_{j1})^3}{(\sum_jn_j^{-1}\var
X_{j1})^{3/2}}
\end{equation}
is a measure of standardised skewness and, under our assumptions, is
bounded as $n\rai$. Result (\ref{eq2.2}) is a version of (\ref
{eq1.1}) in a particular case.

However, in general (\ref{eq2.2}) does not hold in the lattice-valued
case that
we are considering. For example, if $k\geq1$ and the $X_{ji}$s, for all
$i$ and $j$, have a common lattice distribution, then, as was made
clear by Esseen \cite{Ess45}, any expansion of the distribution of $X$
has to include a discontinuous term of size $n^{-1/2}$ (specifically, the
term $n^{-1/2}d_n(x) \phi(x)$ in~(\ref{eq1.2})) that reflects the ``continuity
correction'' needed to approximate the discontinuous distribution of
$T=\{S-E(S)\}/(\var S)^{1/2}$ by a continuous normal distribution.

When exploring this problem, we suppose that the sequence of values of
$n$ is strictly increasing. Further, we assume that
%
%
\begin{equation}
\label{eq2.4} \min_{1\leq j\leq k} \liminf_{n\rai} (n_j/n)>0.
\end{equation}
In Theorem~\ref{the1}, below, we fix both $k$ and the distributions of
$X_{j1}$,
for $1\leq j\leq k$. This means that $e_1,\ldots,e_k$ are fixed too.
However, for each $n$ we consider there to be a potentially new
sequence of values $n_1,\ldots,n_k$. In particular, the ratios
$n_{j_1}/n_{j_2}$ can change considerably from one choice of $n$ to
another, although, in view of (\ref{eq2.4}), $n_{j_1}/n_{j_2}$ is
bounded away
from zero and infinity as $n\rai$.

In the first part of Theorem~\ref{the1}, below, we also impose the following
condition on at least one of the ratios $\rho_{j_1j_2}=(e_{j_2}n_{j_1})/(e_{j_1}n_{j_2})$:
%
%
\begin{equation}
\label{eq2.5} \mbox{for each integer $\ell\geq1$}, \qquad n^{1/2}|\sin(
\ell\rho_{j_1j_2}\pi)|\rai
\end{equation}
as $n\rai$.

\begin{them}\label{the1}
Assume that $E|X_{j1}|^3<\infty$ for $j=1,\ldots,k$; that
$X_{j1}$ is distributed on a lattice $x_j+\nu e_j$, for integers $\nu
$, where $e_j$ is the maximal lattice edge width; and that (\ref{eq2.4})
holds. (\textup{i})~If, for some pair $j_1,j_2$ with $1\leq j_1<j_2\leq k$,
$\rho_{j_1j_2}$ satisfies (\ref{eq2.5}), then the one-term Edgeworth
expansion at (\ref{eq2.2}) holds uniformly in $x$. (\textup{ii}) However, if
$\rho_{j_1j_2}$ equals a fixed rational number (not depending on $n$) for
each pair $j_1,j_2$, and if the points $x_j$ can all be taken equal,
then the expansion at (\ref{eq2.2}) fails to hold because it does not include
an appropriate discontinuous term of size $n^{-1/2}$.
\end{them}

\subsection{\texorpdfstring{Circumstances where (\protect\ref{eq2.5}) holds}
{Circumstances where (2.5) holds}}\label{sec2.2}

If $\rho_0$ is irrational, then $|\sin(\ell\rho_0\pi)|>0$ for all
integers $\ell$. Therefore, (\ref{eq2.5}) holds if $\rho_{j_1j_2}$ converges
to an irrational number as $n\rai$.

However, in many cases (\ref{eq2.5}) holds without the sequence $\rho_{j_1j_2}$ converging. For example, assume for simplicity that the
lattice edge widths $e_j$ are all identical, let $\rho_1$ and $\rho_2$ be two distinct irrational numbers, and let the sequence of values
of the ratio $n_{j_1}/n_{j_2}$ be a sequence of convergents of $\rho_1$ and $\rho_2$, chosen so that an infinite number of convergents
come from each $\rho_j$. (For a definition of convergents of
irrational numbers, see, e.g., Leveque \cite{Lev56}, p. 70.) Then
(\ref{eq2.5}) holds, although the sequence $\rho_{j_1j_2}$ does not converge.

Importantly, (\ref{eq2.5}) also holds in many cases where each $\rho_{j_1j_2}$
is close to a rational number, indeed where each $\rho_{j_1j_2}$
converges to a rational number. For example, we claim that (\ref{eq2.5})
obtains if $\rho_{j_1j_2}=1+\ep_{j_1j_2}$, where $\ep_{j_1j_2}=\ep_{j_1j_2}(n)$,
which\vspace*{1pt} can be either positive or negative, converges to
zero strictly more slowly than $n^{1/2}$:
%
%
\begin{equation}
\label{eq2.6} \ep_{j_1j_2}\ra0, \qquad n^{1/2}|
\ep_{j_1j_2}| \rai.
\end{equation}
In this case, for each fixed, positive integer $\ell$,
\[
\sin(\ell\rho_{j_1j_2}\pi)=\sin(\ell\pi)+\ell\pi\ep_{j_1j_2} \cos(
\ell\pi)+O \bigl(\ep_{j_1j_2}^2 \bigr),
\]
from which it follows that
%
%
\begin{equation}
\label{eq2.7} \bigl|\sin(\ell\rho_{j_1j_2}\pi)\bigr| \sim %
\cases{ \ell
\pi| \ep_{j_1j_2}| & if $\ell$ is an even integer
\cr
1 & if $\ell$ is an
odd integer,}
\end{equation}
where $a_n\sim b_n$ means that the ratio $a_n/b_n$ converges to 1.
Assumption (\ref{eq2.5}) follows from (\ref{eq2.6}) and (\ref{eq2.7}).

A similar argument can be used to prove that if $\rho_{j_1j_2}=\rho_0+\ep_{j_1j_2}$, where $\rho_0$ is a fixed rational number and $\ep_{j_1j_2}=\ep_{j_1j_2}(n)$ satisfies (\ref{eq2.6}), then (\ref
{eq2.5}) is true. (The
case $\rho_0=0$ is excluded by (\ref{eq2.4}).) These examples make it clear
that there is not a great deal of latitude in the assumption, imposed
in part (ii) of Theorem~\ref{the1}, that each $\rho_{j_1j_2}$ should
equal a
fixed rational number. In particular, for (\ref{eq2.5}) to fail it is not
sufficient that each $\rho_{j_1j_2}$ converge to a rational.

\subsection{Refinement of bound on remainder term in Edgeworth expansions}\label{sec2.3}

In Section~\ref{sec2.1}, we showed that, if (\ref{eq2.5}) holds, the
discontinuous term
of size $n^{-1/2}$, in expansions such as (\ref{eq1.2}), is actually
of smaller
order than $n^{-1/2}$. To obtain a more concise bound on the
discontinuous term, we shall investigate in detail cases where one or
more of the ratios $\rho_{j_1j_2}$ converge to an irrational number as
$n$ diverges. However, this treatment requires a definition of the
``type'' of an irrational, and we give that next.

If $x$ is a real number, let $\langle x\rangle$ denote the distance
from $x$
to the nearest integer. (In particular, if $\lfloor x\rfloor$ is the integer
part function, $\langle x\rangle=\min\{x-\lfloor x\rfloor
,1-(x-\lfloor x\rfloor)\}$.)
We say that the irrational number $\rho$ is of type $\eta$ if $\eta$
equals the supremum of all $\ze$ such that $\liminf_{p\rai} p^\ze
\langle p\rho\rangle=0$, where $p\rai$ through integer values. Properties
of convergents of irrational numbers (specifically, Dirichlet's
Theorem) can be used to prove that the type of any given irrational
number always satisfies $\eta\geq1$. It follows from Roth's Theorem
(Roth \cite{Rot55}) that all algebraic irrationals (that is, all
irrational numbers that are roots of polynomials with rational
coefficients) are of minimal type, i.e., $\eta=1$, which is one of the
cases we consider below.

More generally, if a number is chosen randomly, for example as the
value of a random variable having a continuous distribution on the real
line, then with probability 1 it is an irrational of type~1.
Irrationals that are not algebraic are said to be transcendental, and
can have type strictly greater than 1. (However, the transcendental
number $e$ is of type 1.) Known upper bounds to the types of $\pi$,
$\pi^2$ and $\log2$ are 6.61, 4.45 and 2.58, respectively. Liouville
numbers have type $\eta=\infty$. The type of an irrational number is
one less than its irrationality measure (or equivalently, one less than
its approximation exponent or Liouville-Roth constant). We refer the
reader to Ribenboim \cite{Rib00} for more information about types of
irrational numbers.

Next, we introduce notation which helps us to define an approximation
to the discontinuous term, an analogue of $d_n(x)$ in (\ref{eq1.2}), when
$k=2$. (Here, $k$ is as in (\ref{eq2.1}).) Assuming that the lattice,
on which
the distribution of $X_{ji}$ is supported, consists of points $x_j+\nu
e_j$ for integers $\nu$, define $\xi_{jn}=e_j (\si_j
n_j^{1/2})^{-1}\{
(n_j x_j/e_j)-\lfloor n_j x_j/e_j\rfloor\}$ and
%
%
\begin{equation}
\label{eq2.8} \xi_n(x)= \bigl\{x-(c_1
\xi_{1n}+c_2 \xi_{2n}) \bigr\} \frac{\si_1 n_1^{1/2}}{
c_1 e_1},
\end{equation}
where, recalling that $\si_j^2=\var(X_{ji})$, we define $c_j$ for
$j=1$ and 2 by
%
%
\begin{equation}
\label{eq2.9} c_j= \biggl(\frac{n_j^{-1}\si_j^2}{ n_1^{-1}\si_1^2+n_2^{-1}\si_2^2} \biggr)^{ 1/2}.
\end{equation}

Let $\al\in(0,{\frac{1}{2}})$ and partition the set of all integers into
adjacent blocks each comprised of $2 \lfloor n^\al\rfloor+1$ consecutive
integers. Write ${\bar\nu}_\ell$ for the central integer in the
$\ell$th
block, which we denote by $\cN_\ell$ where $-\infty<\ell<\infty$
and $\cN_{\ell+1}$ is immediately to the right of $\cN_\ell$ on the
number line. Given $\nu\in\cN_\ell$, put $\nu_\ell=\nu-{\bar\nu}_\ell$.

Let $c_3=e_2 n_1/\si_1 n_2$ and $c_4=(e_1/\si_2) (n_1/n_2)^{1/2}$,
and note that $c_1,\ldots,c_4$ are strictly positive functions of $n$
and are bounded away from zero and infinity as $n$ diverges. Put $\ga
=\prod_{j=1,2} (e_j/\si_j)$, and, given an integer $r_0\geq1$, define
%
%
\begin{eqnarray}
\label{eq2.10} \phi(u,x)&=&\phi \bigl\{(x/c_1)-c_3 u
\bigr\} \phi(c_4 u), \qquad\phi_r(u,x)=(\partt /\partt
u)^r \phi(u,x),
\nonumber
\\
K_n(x)&=&\ga\sum_{r=0}^{r_0}
\sum_{-\infty<\ell<\infty} \frac
{\phi_r ({\bar\nu}_\ell /n_1^{1/2},x )}{ r! n_1^{r/2}} \sum
_{\nu\in\cN_\ell} \nu_\ell^r \psi \biggl\{
\xi_n(x)-\frac{e_2 n_1}{
e_1 n_2} \nu \biggr\},
\end{eqnarray}
where, as in Section~\ref{sec1.1}, $\psi(x)=\lfloor x\rfloor
-x+{\frac{1}{2}}$.

We claim that the infinite series in the definition of $K_n(x)$ is
absolutely convergent, uniformly in $x$. To appreciate why, note that
%
%
\begin{equation}
\label{eq2.11} \sup_{-\infty<x<\infty} \bigl|\phi_r \bigl({\bar
\nu}_\ell /n_1^{1/2},x \bigr) \bigr| \leq
C_1(r) \phi \bigl(c_4 {\bar\nu}_\ell
/n_1^{1/2} \bigr),
\end{equation}
where, here and below, the notation $C_j(r)$ will denote a constant
depending on $r$ but not on $n$. Using (\ref{eq2.4}) and (\ref
{eq2.11}), we deduce that
%
%
\begin{equation}
\label{eq2.12} \sum_{-\infty<\ell<\infty} \Bigl\{\sup_{-\infty<x<\infty}
\bigl|\phi_r \bigl({\bar\nu}_\ell /n_1^{1/2},x
\bigr) \bigr| \Bigr\} \leq C_2(r) n^{(1/2)-\al}.
\end{equation}
(In more detail, without loss of generality the block $\cN_0$ is
centred at 0, in which case, when bounding the series on the left-hand
side of (\ref{eq2.12}), ${\bar\nu}_\ell$ can be interpreted as $2
\ell n^\al$.
Consequently the left-hand side of (\ref{eq2.12}) is bounded by a constant
multiple of $C_1(r)\int\phi(2 u n^\al/n_1^{1/2}) \,du$, and (\ref
{eq2.12}) follows.)

More simply, since (a) $|\nu_\ell|\leq n^\al$, (b) $|\cN_\ell|\leq
(2 n^\al+1)$ (where we define $|\cN_\ell|=\#\cN_\ell$), and (c)
$|\psi|\leq{\frac{1}{2}}$, then
%
%
\begin{equation}
\label{eq2.13} \sup_{-\infty<x<\infty} \biggl|\sum_{\nu\in\cN_\ell}
\nu_\ell^r \psi \biggl\{\xi_n(x)-
\frac{e_2 n_1}{ e_1 n_2} \nu \biggr\} \biggr| \leq C_3(r) n^{(r+1) \al}.
\end{equation}
Combining (\ref{eq2.10}), (\ref{eq2.12}) and (\ref{eq2.13}), and
replacing each summand on the
right-hand side of (\ref{eq2.10}) by its absolute value, we obtain the bound:
$n^{-1/2}|K_n(x)|\leq C_4(r_0)$, uniformly in~$x$. This inequality
demonstrates the claimed absolute convergence of the series in~(\ref{eq2.10}).

Recall the definition of $S$ at (\ref{eq2.1}), and that $\rho_{j_1j_2}=(e_{j_2}n_{j_1})/(e_{j_1}n_{j_2})$. Part (i) of Theorem~\ref{the2},
below, captures the analogue of the discontinuous term, $d_n(x)$, in a
multisample version of (\ref{eq1.2}), and part (ii) gives conditions under
which the net contribution of that term equals $O(n^{\delta
-(1/2)-(1/2\eta)})$, for all $\delta>0$ when some $\rho_{j_1j_2}$ is
sufficiently close to an irrational number of type $\eta$.

\begin{them}\label{the2}
Assume that $E|X_{j1}|^4<\infty$ for $j=1,\ldots,k$; that
$X_{j1}$ is distributed on a lattice $x_j+\nu e_j$, for integers $\nu
$, where $e_j$ is the maximal lattice edge width; and that (\ref{eq2.4})
holds. Choose $r_0\geq4\al/(1-2\al)$ in (\ref{eq2.10}). (\textup{i}) If
$k=2$ and
$K_n$ is as defined at (\ref{eq2.10}), then
%
%
\begin{equation}
\label{eq2.14}
\hspace*{-5pt}P \biggl\{\frac{S-E(S)}{(\var S)^{1/2}}\leq x \biggr\} =\Phi(x)+n^{-1/2}{
\frac{1}{6}}\be \bigl(1-x^2 \bigr) \phi(x) +(n_1n_2)^{-1/2}K_n(x)
+O \bigl(n^{-1} \bigr),
\end{equation}
uniformly in $x$. (\textup{ii}) If, for some pair $j_1,j_2$ with $1\leq
j_1<j_2\leq k$, the ratio $\rho_{j_1j_2}=(e_{j_2}n_{j_1})/ (e_{j_1}n_{j_2})$ satisfies
%
%
\begin{equation}
\label{eq2.15} |\rho_{j_1j_2}-\rho_0|=O \bigl(n^{-(1/2) \{1+(1/\eta)+\delta\}}
\bigr)
\end{equation}
for some $\delta>0$, where $\rho_0$ is an irrational number of type
$\eta$, then, for each $\delta>0$,
%
%
\begin{equation}
\label{eq2.16} P \biggl\{\frac{S-E(S)}{(\var S)^{1/2}}\leq x \biggr\} =\Phi(x)+n^{-1/2}{
\frac{1}{6}}\be \bigl(1-x^2 \bigr) \phi(x) +O
\bigl(n^{\delta-(1/2)-(1/2\eta)} \bigr),
\end{equation}
uniformly in $x$.
\end{them}

Result (\ref{eq2.16}) is of particular interest in the case $\eta=1$, which
encompasses almost all irrational numbers (with respect to Lebesgue
measure), including all the algebraic irrationals and some
transcendental numbers. When $\eta=1$,
%
%
\begin{equation}
\label{eq2.17} P \biggl\{\frac{S-E(S)}{(\var S)^{1/2}}\leq x \biggr\} =\Phi(x)+n^{-1/2}{
\frac{1}{6}}\be \bigl(1-x^2 \bigr) \phi(x) +O
\bigl(n^{\delta-1} \bigr),
\end{equation}
uniformly in $x$ for each $\delta>0$. Result (\ref{eq2.17}) implies
that the
lattice nature of the distribution of $X_{ji}$ can be ignored, almost up
to terms of second order in Edgeworth expansions, when considering the
impact of latticeness on the accuracy of normal approximations.

\subsection{Practical choice of $n_1$ and $n_2$}\label{sec2.4}

In practice it is not difficult to choose $n_1$ and $n_2$ so that
(\ref{eq2.15}) holds. To see how, assume for simplicity that the
lattice edge
widths $e_1$ and $e_2$ are identical, as they would be if (for example)
$S$ were equal to a sum or difference of estimators of binomial
proportions. If $\rho_0$ is an irrational number then the convergents
$m_1/m_2$ of $\rho_0$ satisfy
%
%
\begin{equation}
\label{eq2.18} \bigl|(m_1/m_2)-\rho_0\bigr|\leq
m_2^{-2}.
\end{equation}
(See e.g. Leveque \cite{Lev56}, equation (29), p. 180.) Therefore, if
$n_1$ and $n_2$ are relatively prime and $n_1/n_2$ is a convergent of
$\rho_0$, then (\ref{eq2.15}), for each $\delta\in(0,3-(1/\eta)]$, follows
from (\ref{eq2.18}). The most difficult case, as far as (\ref
{eq2.15}) is concerned, is
the one where the convergence rate in (\ref{eq2.15}) is fastest, and arises
when $\eta=1$. There we need to ensure that
%
%
\begin{equation}
\label{eq2.19} |\rho_{j_1j_2}-\rho_0|=O \bigl(n^{-1-\delta}
\bigr)
\end{equation}
for some $\delta>0$. Now, (\ref{eq2.19}) holds whenever $n_1/n_2$ is a
convergent of $\rho_0$, and the Khinchin-L\'evy Theorem (see, e.g.,
pp. 82--83 of Einsiedler and Ward \cite{Ein11}) implies that the
convergents are reasonably closely spaced; the numerators and
denominators generally increase by factors of only $\pi^2/(12 \log
2)\approx1.87$. Moreover, there are many ratios $n_1/n_2$ on either
side of convergents for which (\ref{eq2.19}) holds.

The pair $(n_1,n_2)$ can be chosen from tables of, or formulae for,
convergents for commonly arising irrationals of type 1. See, for
example, Griffiths \cite{Gri04} and references therein, and note that
$e$ and any algebraic irrational is of type 1.

\subsection{\texorpdfstring{Alternative formula for $K_n$, and derivation of (\protect\ref{eq2.16})
from (\protect\ref{eq2.14}) when~$\eta=1$}
{Alternative formula for $K_n$, and derivation of (2.16) from (2.14) when eta=1}}\label{sec2.5}

Part (i) of Theorem~\ref{the2} can be stated for a version of $K_n(x)$ simpler
than that at (\ref{eq2.10}):
%
%
\begin{equation}
\label{eq2.20} K_n(x)=\ga\sum_\nu\phi
\biggl(\frac{x}{ c_1} -\frac{e_2 n_1^{1/2}}{\si_1 n_2} \nu \biggr) \phi \bigl
\{e_2 \bigl(\si_2 n_2^{1/2}
\bigr)^{-1}\nu \bigr\} \psi \biggl\{\xi_n(x)-
\frac{e_2 n_1}{ e_1 n_2} \nu \biggr\}.
\end{equation}
Indeed, the $K_n(x)$ at (\ref{eq2.20}) is just $\ga I_4(x)$, where
$I_4(x)$ is
as defined at (\ref{eq4.17}) in the proof of Theorem~\ref{the1}, and
in fact that
formula provides a convenient point of access to a proof of (\ref
{eq2.14}) with
$K_n(x)$ as at (\ref{eq2.20}). However, in the case $\eta>1$ it is not
straightforward to pass from (\ref{eq2.20}) to (\ref{eq2.16}), and
that is why we used
the definition of $K_n(x)$ at (\ref{eq2.10}).

To appreciate that (\ref{eq2.16}) follows from (\ref{eq2.20}) when
$\eta=1$, note that
the definition of $K_n(x)$ at (\ref{eq2.20}) is equivalent to:
%
%
\begin{equation}
\label{eq2.21} \ga^{-1}K_n(x)=\sum
_\nu\Psi(x,\nu) \psi \biggl\{\xi_n(x)-
\frac{e_2 n_1}{ e_1 n_2} \nu \biggr\},
\end{equation}
where
\[
\Psi(x,\nu)=\phi \biggl(\frac{x}{ c_1} -\frac{e_2 n_1^{1/2}}{\si_1 n_2} \nu \biggr) \phi
\bigl\{e_2 \bigl(\si_2 n_2^{1/2}
\bigr)^{-1}\nu \bigr\}.
\]
If (\ref{eq2.15}) holds with $\eta=1$ then a standard argument for bounding
discrepancies of sequences (see p. 123 of Kuipers and Niederreiter
\cite{Kui74}) can be used to prove that for all $\delta>0$,
%
%
\begin{equation}
\label{eq2.22} \sup_{-\infty<z<\infty} \Biggl|\sum_{\nu=1}^N
\psi \biggl(z-\frac{e_2 n_1}{ e_1 n_2} \nu \biggr) \Biggr| =O \bigl(N^\delta \bigr).
\end{equation}
Note too that
%
%
\begin{equation}
\label{eq2.23} \sup_{\nu\geq1} \sup_{-\infty<x<\infty} \bigl|\Psi(x,\nu+1)-\Psi(x,
\nu)\bigr|\leq C n^{-1/2}.
\end{equation}
Taking $a_\nu=\Psi(x,\nu)$ and $b_\nu=\psi\{\xi_n(x)-(e_2 n_1/e_1
n_2) \nu\}$, and employing Abel's method of summation, we can write:
\[
\sum_{\nu=1}^N a_\nu
b_\nu =a_N \sum_{\nu=1}^N
b_\nu-\sum_{\nu=1}^{N-1}
(a_{\nu+1}-a_\nu) \sum_{j=1}^\nu
b_j,
\]
which in company with (\ref{eq2.22}) and (\ref{eq2.23}) allows us to
prove that,
provided $N=O(n^C)$ for some $C>0$,
%
%
\begin{equation}
\label{eq2.24} \sup_{-\infty<x<\infty} \Biggl|\sum_{\nu=1}^N
\Psi(x,\nu) \psi \biggl\{\xi_n(x)-\frac{e_2 n_1}{ e_1 n_2} \nu \biggr\}\Biggr |
=O \bigl(N^\delta \bigr)
\end{equation}
for all $\delta>0$. More simply, if $N\geq n^2$ then
%
%
\begin{eqnarray}
\label{eq2.25} &&\sup_{-\infty<x<\infty}\Biggl |\sum_{\nu=N+1}^\infty
\Psi(x,\nu) \psi \biggl\{\xi_n(x)-\frac{e_2 n_1}{ e_1 n_2} \nu \biggr\} \Biggr|
\nonumber
\\
&&\hspace*{157pt}\leq\frac{1}{2} \sup_{-\infty<x<\infty} \sum
_{\nu=N+1}^\infty \Psi(x,\nu) =O(1).
\end{eqnarray}
Combining (\ref{eq2.24}) and (\ref{eq2.25}), using a similar argument
to treat series
where $\nu\leq0$, and noting the definition of $K_n(x)$ at (\ref
{eq2.21}), we
deduce that $\sup_x |K_n(x)|=O(n^\delta)$ for all $\delta>0$. In the
case $k=2$, and for $\eta=1$, this gives (\ref{eq2.16}) as a
corollary of (\ref{eq2.14}).

\subsection{\texorpdfstring{Derivation of (\protect\ref{eq2.16}) from (\protect\ref{eq2.14}) when $\eta\geq1$}
{Derivation of (2.16) from (2.14) when eta>=1}}\label{sec2.6}

Our proof of (\ref{eq2.16}), in Section~\ref{sec4.2}, will proceed by deriving
implicitly a version of (\ref{eq2.14}) in the case $k\geq2$, and
showing that,
if (\ref{eq2.15}) holds, then that version of (\ref{eq2.14}) entails
(\ref{eq2.16}). The
relative complexity of a form of (\ref{eq2.14}) for general $k$
discouraged us
from including it in Theorem~\ref{the2}, but it is nevertheless
instructive to
show how, when $k=2$, one can obtain (\ref{eq2.16}) from (\ref
{eq2.14}). We outline the
proof below, highlighting the properties of irrational numbers,
particularly the differences between the case of irrationals of type
$\eta=1$ and the case of those of larger type, that determine the
bound for the remainder in (\ref{eq2.16}).

Note that if $q$ is a polynomial function then, applying Koksma's
inequality (see, e.g., Theorem~5.1, p. 143 of \cite{Kui74}) and the
Erd\H os-Tur\'an inequality (see, e.g., formula (2.42), p. 114 of \cite
{Kui74}), it can be shown that
%
%
\begin{equation}
\label{eq2.26} \chi(N,q,\tau)\equiv\sup_{-\infty<z<\infty} \Biggl|\sum
_{i=1}^N q(i/N) \psi(z-\tau i)\Biggr | \leq
C_1(q) \Biggl\{\frac{N}{ m} +\sum
_{\ell=1}^m \frac{1}{\ell|\sin(\ell\tau\pi)|} \Biggr\} ,
\end{equation}
for all integers $m\geq1$. Here $\tau>0$ is permitted to vary with
$N$, and the constant $C_1(q)$ depends on the degree and the
coefficients of $q$ but not on the positive integer $N$ or on $m$, $z$
or $\tau$.

We shall take $\tau=\rho_{12}$, a function of $n$, in which case, since
\[
\bigl|\bigl|\sin(\ell\tau\pi)\bigr|-\bigl|\sin(\ell\rho_0\pi)\bigr| \bigr|\leq\ell \pi|
\rho_{12}-\rho_0|,
\]
we have:
%
%
\begin{equation}
\label{eq2.27} \bigl|\sin(\ell\tau\pi)\bigr|+\ell\pi|\rho_{12}-
\rho_0| \geq\bigl|\sin(\ell\rho_0\pi)\bigr| =\sin \bigl(\pi\langle
\ell\rho_0\rangle \bigr) \geq2 \langle\ell\rho_0\rangle
\geq2 C_2(\delta) \ell^{-\eta-\delta}
\end{equation}
for any given $\delta>0$ and all $\ell\geq1$, where $C_2(\delta)>0$
depends on $\delta$ but not on $\ell$, the last inequality in (\ref{eq2.27})
follows from the assumption that $\rho_0$ is of type $\eta$, and the
second-last inequality comes from the fact that $0\leq\langle x\rangle
\leq
{\frac{1}{2}}$ for all real numbers $x$, and $\sin(\pi x)\geq2 x$ whenever
$0\leq x\leq{\frac{1}{2}}$. If
%
%
\begin{equation}
\label{eq2.28} |\rho_{12}-\rho_0|\leq C_2(
\delta) \pi^{-1}m^{-(1+\eta+\delta)}
\end{equation}
then it follows from (\ref{eq2.27}) that $|\sin(\ell\tau\pi)|\geq
\langle
\ell
\rho_0\rangle$ for $1\leq\ell\leq m$, and so
%
%
\begin{equation}
\label{eq2.29} \sum_{\ell=1}^m
\frac{1}{\ell|\sin(\ell\tau\pi)|} \leq\sum_{\ell=1}^m
\frac{1}{\ell\langle\ell\rho_0\rangle} .
\end{equation}

A standard argument for bounding the discrepancy of a sequence (see,
e.g., p. 123 of~\cite{Kui74}) can be used to show that, since $\rho_0$ is an irrational number of type $\eta$,
%
%
\begingroup
\abovedisplayskip=6.5pt
\belowdisplayskip=6.5pt
\begin{equation}
\label{eq2.30} \sum_{\ell=1}^m
\frac{1}{\ell\langle\ell\rho_0\rangle} =O \bigl(m^{\eta-1+\delta} \bigr)
\end{equation}
for all $\delta>0$. Therefore, provided that (\ref{eq2.28}) holds, we can
deduce from (\ref{eq2.26}) and (\ref{eq2.29}) that
%
%
\begin{equation}
\label{eq2.31} \chi(N,q,\rho_{12}) \leq C_3(q)
m^{-1} \bigl(N+m^{\eta+\delta} \bigr) =O \bigl(N^{1-(1/\eta)+\delta_1}
\bigr),
\end{equation}
where $\delta_1>0$ and the inequality holds for all $m$ and the
identity is true if $m/N^{1/\eta}$ is bounded away from zero and
infinity as $N\rai$. When $m$ has the latter property, (\ref{eq2.28}) is
satisfied, for all sufficiently large $N$, if
%
%
\begin{equation}
\label{eq2.32} |\rho_{12}-\rho_0|=O \bigl(N^{-\{1+(1/\eta)+\delta_2\}}
\bigr)
\end{equation}
for some $\delta_2>\delta_1$.

Note that it is at (\ref{eq2.30}) that the type, $\eta$, of the irrational
number $\rho_0$ enters into consideration. In the case $\eta=1$ the
exponent $\delta$ in (\ref{eq2.30}) could not be removed or reduced, perhaps
by replacing the implicit factor $m^\delta$ in (\ref{eq2.30}) by
$(\log m)^C$
for some $C>0$, without an analogous strengthening of Roth's Theorem.
Formula (\ref{eq2.30}) also marks the step at which it becomes
apparent that a
poorer bound will be obtained for an irrational number of type 1,
relative to one of type $\eta>1$.

Applying the bound (\ref{eq2.31}), for several versions of the
polynomial $q$,
in the case $N=2 \lfloor n^\al\rfloor+1$, we deduce that
%
%
\begin{equation}
\label{eq2.33} n^{-r/2} \sup_{-\infty<z<\infty} \biggl|\sum
_{\nu\in\cN_\ell} \nu_\ell^r \psi \biggl(z-
\frac{e_2 n_1}{ e_1 n_2} \nu \biggr) \biggr| =O \bigl(n^{r \{\al-(1/2)\}+\al\{1-(1/\eta)+\delta_1\}} \bigr) ,
\end{equation}
provided that (\ref{eq2.32}) holds, i.e., $|\rho_{12}-\rho_0|=O(n^{-\al\{
1+(1/\eta)+\delta_2\}})$. Now, the only constraint on $\al$ is
$0<\al<{\frac{1}{2}}$, and so we can choose $\al$ as close to
${\frac{1}{2}}$, but
less than ${\frac{1}{2}}$, as we desire. In particular, if $\delta_3>0$ is
given, and we choose $\al={\frac{1}{2}}-\delta_4$ where $\delta_4>0$ is
sufficiently small, then by taking $\delta_1$ in (\ref{eq2.33}) to be
small we obtain:
%
%
\begin{equation}
\label{eq2.34} \max_{1\leq r\leq r_0} n^{-r/2} \sup_{-\infty<z<\infty} \biggl|\sum
_{\nu\in\cN_\ell} \nu_\ell^r \psi
\biggl(z-\frac{e_2 n_1}{ e_1 n_2} \nu \biggr) \biggr| =O \bigl(n^{(1/2) \{1-(1/\eta)\}+\delta_3} \bigr),
\end{equation}
provided that
%
%
\begin{equation}
\label{eq2.35} |\rho_{12}-\rho_0|=O \bigl(n^{-(1/2) \{1+(1/\eta)+\delta_5\}}
\bigr) ,
\end{equation}
where $\delta_5>0$ can be made as small as we like simply by choosing
$\delta_4$ small. Now, (\ref{eq2.35}) follows from (\ref{eq2.15}).
It therefore
follows from (\ref{eq2.34}), and the definition of $K_n(x)$ at (\ref
{eq2.10}), that if
(\ref{eq2.15}) holds for some $\delta>0$ then
%
%
\begin{equation}
\label{eq2.36} \sup_{-\infty<x<\infty} (n_1n_2)^{-1/2}\bigl|K_n(x)\bigr|
=O \bigl(n^{\delta-(1/2) \{1+(1/\eta)\}} \bigr)
\end{equation}
\endgroup
for all $\delta>0$. Results (\ref{eq2.14}) and (\ref{eq2.36}) imply
(\ref{eq2.16}), as had to
be shown.\vadjust{\eject}

\subsection{Expansions relating to the bootstrap}\label{sec2.7}

In this section we show that, despite the potential for problems
arising from discreteness, the bootstrap (including the double
bootstrap) applied to inference based on the distribution of $\{S-E(S)\}
/(\var S)^{1/2}$, generally (when (\ref{eq2.15}) holds and $\rho_0$
is of type
1) produces confidence regions and hypothesis tests with the same
orders of magnitude of coverage or level accuracy, up to terms of size
$n^{\delta-1}$ for all $\delta>0$, as it would in the case of smooth
sampling distributions. This result is of practical importance, since
standard percentile bootstrap methods applied to lattice distributions
are frustrated by the effects of discontinuities; see, e.g., Singh
\cite{Sin81} and Hall \cite{Hal89}.

For brevity, when establishing this property we treat only the context
of Theorem~\ref{the2}. We begin by stating an analogue of (\ref
{eq2.14}) there, valid
when $k=2$. The arguments used to prove part (i) of Theorem~\ref{the2}
can be
employed to show that
%
%
\begin{eqnarray}
\label{eq2.37} &&P \biggl[\frac{S^*-E(S^*\mid\cX)}{\{\var(S\mid\cX)\}^{1/2}}\leq x \Bigm|\cX \biggr] =
\Phi(x)+n^{-1/2}{\frac{1}{6}} {\hat\be} \bigl(1-x^2
\bigr) \phi (x)
\nonumber
\\
&&\hspace*{170pt} {}+(n_1n_2)^{-1/2}{\widehat
K}_n(x)+n^{-1}\De_1(x),
\end{eqnarray}
where, analogously to the definitions in Section~\ref{sec2.1},
$S^*=\sum_j{\bar X}_j^*$; ${\bar X}_j^*=n_j^{-1}\sum_iX_{ji}^*$ and $X_{j1}^*,\ldots
,X_{jn_j}^*$ are drawn by sampling randomly, with replacement, from
$\cX_j=(X_{j1},\ldots,X_{jn_j})$; $\cX=(\cX_1,\ldots,\cX_k)$;
\[
{\hat\be}=\frac{n^{1/2}E[\{S^*-E(S^*\mid\cX)\}^3\mid\cX]}{\{
\var(S^*
\mid\cX)\}^{3/2}} ;
\]
using (\ref{eq2.10}) or (\ref{eq2.20}), respectively, as the model
for $K_n(x)$,
\begin{eqnarray*}
{\widehat K}_n(x)&=&{\hat\ga}\sum_{r=0}^{r_0}
\sum_{-\infty<\ell
<\infty} \frac{{\hat\phi}_r ({\bar\nu}_\ell /n_1^{1/2},x )}{
r! n_1^{r/2}} \sum
_{\nu\in\cN_\ell} \nu_\ell^r \psi \biggl\{{\hat
\xi}_n(x)-\frac{e_2 n_1}{ e_1 n_2} \nu \biggr\} ,
\\
{\widehat K}_n(x)&=&{\hat\ga}\sum_\nu
\phi \biggl(\frac{x}{{\hat c}_1} -\frac{e_2 n_1^{1/2}}{{\hat\si}_1 n_2} \nu \biggr) \phi \bigl
\{e_2 \bigl({\hat\si}_2 n_2^{1/2}
\bigr)^{-1}\nu \bigr\} \psi \biggl\{{\hat\xi}_n(x)-
\frac{e_2 n_1}{ e_1 n_2} \nu \biggr\},
\end{eqnarray*}
where ${\hat\ga}=\prod_{j=1,2}(e_j/{\hat\si}_j)$,
${\hat\phi}_r(u,x)=(\partt/\partt u)^r {\hat\phi}(u,x)$,
${\hat\phi}(u,x)=\phi\{(x/{\hat c}_1)-{\hat c}_3 u\}\times\phi
({\hat c}_4 u)$,
\[
{\hat c}_j= \biggl(\frac{n_j^{-1}{\hat\si}_j^2}{ n_1^{-1}{\hat\si
}_1^2+n_2^{-1}{\hat\si}_2^2} \biggr)^{ 1/2}
\]
for $j=1$ and 2, ${\hat c}_3=e_2 n_1/{\hat\si}_1 n_2$, ${\hat
c}_4=(e_1/{\hat\si}_2)
(n_1/n_2)^{1/2}$, and ${\hat\xi}_n(x)$ is defined using the empirical
analogue of (\ref{eq2.8}); and, for $C_1>0$ sufficiently large and for
some $C_2>0$,
%
%
\begin{equation}
\label{eq2.38} P \Bigl\{\sup_{-\infty<x<\infty} \bigl|\De_1(x)\bigr|>C_1
n^{-1} \Bigr\}=O \bigl(n^{-C_2} \bigr).
\end{equation}
The assumptions needed for (\ref{eq2.37}) are those imposed for part
(i) of
Theorem~\ref{the2}. The size of $C_2$ in (\ref{eq2.38}) depends to
some extent on the
distributions of $X_{1i}$ and $X_{2i}$ (recall that at this point we
are assuming that $k=2$), but for distributions such as the Bernoulli
or Poisson, which have all moments finite, $C_2$ can be taken
arbitrarily large if $C_1$ is sufficiently large. The connection to
moments here arises because the $O(n^{-C_2})$ bound in (\ref{eq2.38}) is
derived using a method related to Markov's inequality, which can be
applied at a higher order if more moments are finite.

The methods used in Sections~\ref{sec2.5} and~\ref{sec2.6} to derive
uniform bounds to
$K_n$ can also be employed to bound ${\widehat K}_n$, giving
%
%
\begin{equation}
\label{eq2.39} P \Bigl\{\sup_{-\infty<x<\infty} \bigl|{\widehat K}_n(x)
\bigr|>n^{\delta
+(1/2)\{1-(1/\eta)\}} \Bigr\}=O \bigl(n^{-C_3} \bigr)
\end{equation}
for all $\delta>0$ and some $C_3>0$, provided that (\ref{eq2.15})
holds. In
(\ref{eq2.39}), $\eta$ denotes the type of the irrational number
$\rho_0$
appearing in (\ref{eq2.15}), and for sampling distributions such as the
Bernoulli or Poisson (with all moments finite), $C_3$ can be taken
arbitrarily large. Therefore, treating the case of irrationals of type
1, we deduce from (\ref{eq2.37})--(\ref{eq2.39}) that
%
%
\begin{equation}
\label{eq2.40} P \biggl[\frac{S^*-E(S^*\mid\cX)}{\{\var(S^*\mid\cX)\}^{1/2}
}\leq x \Bigm|\cX \biggr] =
\Phi(x)+n^{-1/2}{\frac{1}{6}} {\hat\be} \bigl(1-x^2
\bigr) \phi (x)+n^{-1}\De_2(x) ,
\end{equation}
where
\[
P \Bigl\{\sup_{-\infty<x<\infty} \bigl|\De_2(x)\bigr|>n^\delta \Bigr\}=O
\bigl(n^{-C_2} \bigr) \qquad\mbox{for all $C_2,\delta>0$}.
\]

A similar argument, employing the methods introduced in Section~\ref{sec2.6},
can be used to prove that (\ref{eq2.40}) continues to hold if $k\geq2$,
provided that the assumptions imposed in part (ii) of Theorem~\ref
{the2} hold.
Therefore, the properties stated in the first paragraph of this section hold.

\section{Numerical properties}\label{sec3}

Throughout this section, we take $k=2$ and let $X_{ji}$ be a Bernoulli
random variable satisfying $P(X_{ji}=0)=1-P(X_{ji}=1)=p_j$ for $j=1,2$,
where $p_1=0.4$ and $p_2=0.6$. Thus, $\rho_{12}=e_2 n_1/(e_1
n_2)=n_1/n_2$, where $n_1$ and $n_2$ are the two sample sizes. We take
$n_2$ to be the integer nearest to $\rho_0 n_1$, and vary $n_1$
between 10 and 80; $n_1$ is plotted on the horizontal axes of each of
our graphs. The probability
%
%
\begin{equation}
\label{eq3.1} P(x)=P \bigl[ \bigl\{S-E(S) \bigr\}/(\var S)^{1/2}\leq x
\bigr]
\end{equation}
was approximated by averaging over the results of $10^5$ Monte Carlo
simulations.

%
\begin{figure}

\includegraphics{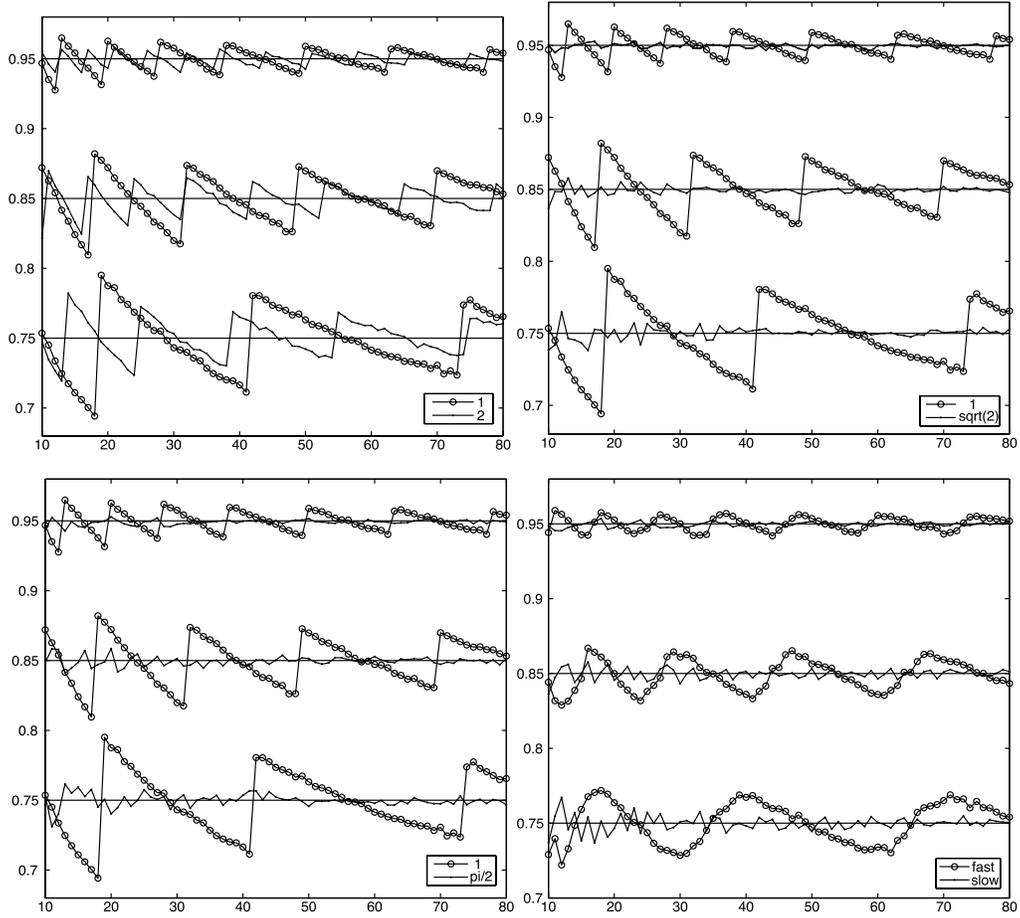}

\caption{Plots of $P(x)$ against $n_1$. Plots are given for $x = \Phi^{-1}(\alpha)=z_\alpha$ and $\alpha=$ 0.95, 0.85, and 0.75, and for
$n_2$ equal to the nearest integer to $\rho_0 n_1$, with $\rho_0=$ 1
or 2 (top left), $\rho_0$= 1 or $2^{1/2}$ (top right), $\rho_0$= 1 or
$\pi$/2 (bottom left) and $\rho_0$ converges to 1 rapidly or slowly
(bottom right; see text for details).}\label{fig1}
\end{figure}

%
\begin{figure}

\includegraphics{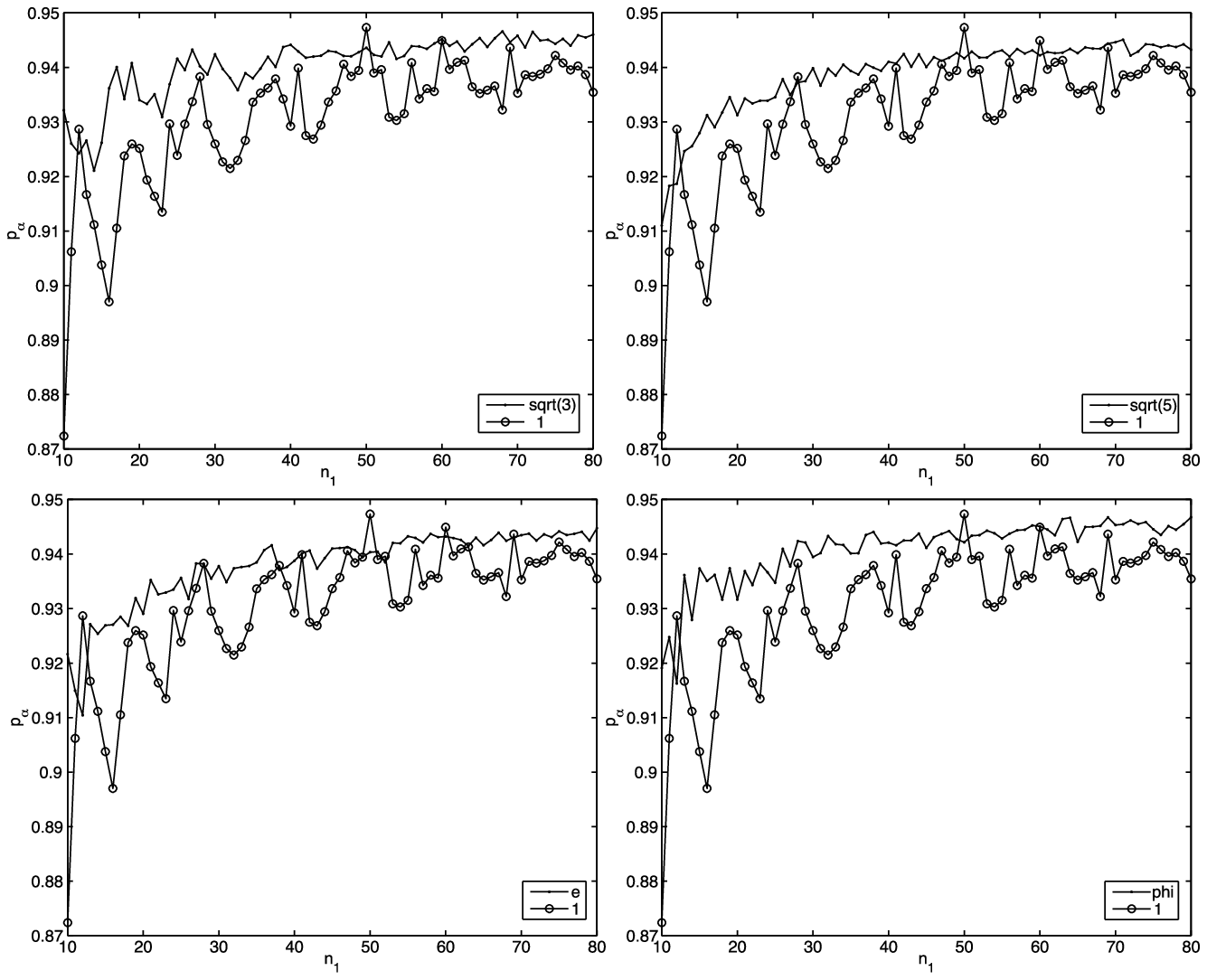}

\caption{Plots of estimates of $P\{E(S)\in\cI_\al\}$, against
$n_1$; see text for details. Each panel shows the case $\rho_0$= 1 and
also, in respective panels, the cases $\rho_0= 3^{1/2}, \rho_0=
5^{1/2}, \rho_0= e$ and $\rho_0= (1 + 5^{1/2})/2$. Throughout, $x
=\Phi^{-1}(\alpha)$ where $\alpha= 0.95$.}\label{fig2}\vspace*{-5pt}
\end{figure}

To illustrate the influence of $\rho_{12}$ on the oscillatory
behaviour of $P(x)$, and in particular on the accuracy of the normal
approximation, each panel in Figure~\ref{fig1} plots $P(x)$ against
$n_1$ for
$x=\Phi^{-1}(\al)=z_\al$ and $\al=0.95$, 0.85 and $0.75$. The top
left panel of Figure~\ref{fig1} shows results for $\rho_0=1$
(indicated by the
lines with circles) and $\rho_0=2$ (lines with dots), and it is clear
that in both cases there is significant oscillatory behaviour, arising
principally from the term in $K_n(x)$ in (\ref{eq2.14}). The top right
panel of
Figure~\ref{fig1} shows that these oscillations decline markedly, and the
accuracy of the normal approximation improves considerably, if $\rho_0=2^{1/2}$. This property reflects the results reported in
Section~\ref{sec2}.

Of course, $\rho_0=2^{1/2}$ is an algebraic irrational. The bottom
left panel of Figure~\ref{fig1} shows that broadly similar values of $P(x)$,
although with somewhat more oscillation (reflecting the relatively low
upper bounds given in Theorem~\ref{the1}), are obtained for $\rho_0=\pi/2$, a
transcendental irrational whose type is bounded above by 6.61. The
bottom right panel of Figure~\ref{fig1} addresses one of the results
reported in
Section~\ref{sec2.2}, specifically that there may be less oscillatory behaviour
when $\rho_{12}$ converges slowly to a rational number than when it
converges quickly. We consider the cases $n_2=n_1+[n_1^{1/5}]$ and
$n_2=n_1+[n_1^{3/5}]$, where $[x]$ denotes the integer nearest to $x$.
In the first case, $\rho_{12}$ converges relatively quickly to 1, and
in the second case the convergence is relatively slow.
Figure~\ref{fig1}
demonstrates that, as anticipated, the oscillatory behaviour is less
pronounced, and the normal approximation better, in the ``slow'' case.

Finally, Figure~\ref{fig2} shows that broadly similar results are
obtained for
coverage probabilities of percentile bootstrap confidence intervals for
$E(S)$. Let $s_\al$ denote the $\al$-level quantile of the
distribution of $S-E(S)$, and let ${\hat s}_\al$, the parametric bootstrap
estimator of $s_\al$, be the $\al$-level quantile of the distribution
of $S^*-S$ given $\cX$, i.e. ${\hat s}_\al=\inf\{s:P(S^*-S\leq
s\mid
\cX)\geq\al\}$.
A~naive $\al$-level one-sided percentile-bootstrap confidence interval
for $E(S)$, with nominal coverage probability $1-\al$, is given by
%
%
\begin{equation}
\label{eq3.3} \cI_\al=(-\infty,S-{\hat s}_\al].\vadjust{\eject}
\end{equation}
In the figure, we give plots of estimates of the coverage probability
$P\{E(S)\in\cI_\al\}$ of $\cI_\al$ against $n_1$, estimated using
$10^5$ Monte-Carlo simulations, for $\al=0.95$. We used $B=9999$
simulations in each bootstrap step. Each panel depicts the case $\rho_0=1$, and successive panels also give results when $\rho_0=3^{1/2}$,
$5^{1/2}$, $e$ and $\phi=(1+5^{1/2})/2$, respectively. Each of these
values of $\rho_0$ is an irrational of type 1, and in each instance
the oscillations are markedly less, and the normal approximation
markedly improved, relative to the case $\rho_0=1$.

%


%

\section{Proofs}\label{sec4}

\subsection{\texorpdfstring{Proof of Theorem~\protect\ref{the1}}{Proof of Theorem 1}}\label{sec4.1}

\subsubsection{\texorpdfstring{Proof of part (i) of Theorem~\protect\ref{the1}}
{Proof of part (i) of Theorem 1}}\label{sec4.1.1}

Here we show that if (\ref{eq2.5}) holds for some $\rho_{j_1j_2}$, where
$j_1\neq j_2$, then (\ref{eq2.2}) obtains. Some of the asymptotic
expansions in
our argument are taken a little further than is necessary for (\ref{eq2.2});
the extra detail will be used in the proof of Theorem~\ref{the2}.

\textit{Step 1: Proof that it is sufficient to consider the case $k=2$.}
Without loss of generality, (\ref{eq2.5}) holds for $\rho_{12}$, and
in this
case we write $S-E(S)=S_1+S_2$, where $S_1=(1-E) ({\bar X}_1+{\bar
X}_2)$ and
$S_2=(1-E) ({\bar X}_3+\cdots+{\bar X}_k)$, where $E$ denotes the expectation
operator. Recall that $S_2$ is independent of ${\bar X}_1$ and ${\bar X}_2$.
Suppose we can prove that, analogously to (\ref{eq2.2}),
%
%
\begin{equation}
\label{eq4.1} P \biggl\{\frac{S_1}{(\var S_1)^{1/2}}\leq x \biggr\} =\Phi(x)+n^{-1/2}{
\frac{1}{6}}\be_1 \bigl(1-x^2 \bigr) \phi(x)+o
\bigl(n^{-1/2} \bigr),
\end{equation}
uniformly in $x$, where, reflecting (\ref{eq2.3}),
\[
\be_1=\be_1(n)=\frac{n^{1/2}E(S_1^3)}{(\var S_1)^{3/2}}.
\]
If we prove that (\ref{eq2.2}), in the case of general $k$, follows
from (\ref{eq4.1}),
we shall have shown that it is sufficient to derive Theorem~\ref{the1}
the case $k=2$.

Since $P(S\leq x)=E\{P(S_1\leq x-S_2\mid S_2)\}$ then we can deduce from
(\ref{eq4.1}) that
%
%
\begin{eqnarray}
\label{eq4.2} &&P(S\leq x)=E \biggl(\Phi \biggl\{\frac{x-S_2}{(\var S_1)^{1/2}} \biggr\} +
\frac{\be_1}{6 n^{1/2}} \biggl[1- \biggl\{\frac{x-S_2}{(\var S_1)^{1/2}} \biggr\}^{ 2}
\biggr]
\nonumber
\\
&&\hspace*{181pt} {}\times\phi \biggl\{\frac{x-S_2}{(\var
S_1)^{1/2}} \biggr\} \biggr)+o
\bigl(n^{-1/2} \bigr),
\end{eqnarray}
uniformly in $x$. Let $R=S_2/(\var S_1)^{1/2}$, and put $\tau_1^2=\var
(R)$, which is bounded away from zero and infinity as $n\rai$. It is
straightforward to prove that, if $N$ denotes a normally distributed
random variable with the same mean (i.e., zero mean) and variance as
$S_2$, then
%
%
\begin{eqnarray}
\label{eq4.3}&&E \biggl( \biggl[1- \biggl\{\frac{x-S_2}{(\var
S_1)^{1/2}} \biggr
\}^{
2} \biggr] \phi \biggl\{\frac{x-S_2}{(\var S_1)^{1/2}} \biggr\} \biggr)
\nonumber
\\
&&\quad=\int \biggl[1- \biggl\{\frac{x}{(\var S_1)^{1/2}}-t \biggr\}^{
2}
\biggr] \phi \biggl\{\frac{x}{(\var S_1)^{1/2}}-t \biggr\} \,dP(R\leq t)
\\
\label{eq4.4}&&\quad=\int \biggl[1- \biggl\{\frac{x}{(\var
S_1)^{1/2}}-t \biggr
\}^{ 2} \biggr] \phi \biggl\{\frac{x}{(\var S_1)^{1/2}}-t \biggr\}
\frac{1}{\tau_1} \phi \biggl(\frac{t}{\tau_1} \biggr)\, dt +O
\bigl(n^{-1/2} \bigr)
\\
\label{eq4.5}&&\quad=E \biggl( \biggl[1- \biggl\{\frac{x-N}{(\var
S_1)^{1/2}} \biggr
\}^{
2} \biggr] \phi \biggl\{\frac{x-N}{(\var S_1)^{1/2}} \biggr\} \biggr) +O
\bigl(n^{-1/2} \bigr),
\end{eqnarray}
uniformly in $x$. The passage from (\ref{eq4.3}) to (\ref{eq4.4}) can
be accomplished
by integrating by parts in (\ref{eq4.3}), then using an Edgeworth
expansion of
the distribution of $R$, then separating out the term in $n^{-1/2}$ in
that expansion, and finally, undoing the integration by parts as it
applies to the leading term in the Edgeworth expansion.

Let $\tau_2^2=\var S_2$ and $\be_2=\be_2(n)=n^{1/2}E(S_2^3)/\tau_2^3$. If
\[
\Phi(r/\tau_2)+n^{-1/2}{\tfrac{1}{6}}\be_2
\bigl\{1-(r/\tau_2)^2 \bigr\} \phi (r/\tau_2)
\]
represents the two-term Edgeworth approximation to $P(S_2\leq r)$ that
would be employed if the distribution of $S_2$ were continuous, then it
can be proved that, uniformly in $x$,
%
%
\begin{eqnarray}
\label{eq4.6} E \biggl[\Phi \biggl\{\frac{x-S_2}{(\var S_1)^{1/2}} \biggr\} \biggr] &=&\int
\Phi \biggl\{\frac{x-r}{(\var S_1)^{1/2}} \biggr\} d_r \biggl\{\Phi(x/
\tau_2)
\nonumber
\\
&&{} +n^{-1/2}{\frac{1}{6}}\be_2 \bigl\{1-(r/
\tau_2)^2 \bigr\} \phi (r/\tau_2) \biggr\}
\nonumber
\\
&&{} + %
\cases{ o \bigl(n^{-1/2} \bigr)& if $
\max_j E|X_{j1}|^3<\infty$
\cr
O
\bigl(n^{-1} \bigr)& if $\max_j E|X_{j1}|^4<
\infty$. } %
\end{eqnarray}
To derive (\ref{eq4.6}), first integrate by parts on the left-hand side,
writing it as
%
%
\begin{eqnarray}
\label{eq4.7} &&\frac{1}{(\var S_1)^{1/2}} \int\phi \biggl\{\frac{x-r}{(\var S_1)^{1/2}} \biggr\}
P(S_2\leq r) \,dr
\nonumber
\\
&&\hspace*{125pt}=\int\phi \biggl\{\frac{x}{(\var S_1)^{1/2}}-t \biggr\} P(R\leq t) \,dt .
\end{eqnarray}
Next, write down an Edgeworth expansion, (E) say, for the joint
distribution of ${\bar X}_3,\ldots, {\bar X}_k$, up to terms of $o(n^{-1/2})$
when $\max_j E|X_{j1}|^3<\infty$ and $O(n^{-1})$ when $\max_j
E|X_{j1}|^4<\infty$. The expansion will include the conventional
discontinuous terms of size $n^{-1/2}$. Use (E) to the derive
discontinuous term $n^{-1/2}D$, say, up to a remainder of smaller order
$n^{-1/2}$, in an Edgeworth expansion of the distribution of $R$. Since
the function $\phi$ is smooth, the impact of $n^{-1/2}D$ on the
right-hand side of (\ref{eq4.7}) equals $O(n^{-1})$, this being
obtained by
multiplying together the factor $n^{-1/2}$ and another term of order
$n^{-1/2}$ that results from integrating $D$ against a smooth function.
Therefore, (\ref{eq4.6}) holds.

Combining (\ref{eq4.2}), (\ref{eq4.5}) and (\ref{eq4.6}), we deduce that
%
%
\begin{eqnarray}
\label{eq4.8} P(S\leq x)&=&\int \biggl(\Phi \biggl\{\frac{x-r}{(\var
S_1)^{1/2}} \biggr\}
+ \frac{\be_1}{6 n^{1/2}} \biggl[1- \biggl\{\frac{x-r}{(\var S_1)^{1/2}} \biggr\}^{ 2}
\biggr] \phi \biggl\{\frac{x-r}{(\var S_1)^{1/2}} \biggr\} \biggr)
\nonumber
\\
&&{} \times d_r \biggl\{\Phi(x/\tau_2)+n^{-1/2}{
\frac{1}{6}}\be_2 \bigl\{ 1-(r/\tau_2)^2
\bigr\} \phi(r/\tau_2) \biggr\} +o \bigl(n^{-1/2} \bigr),
\end{eqnarray}
uniformly in $x$. Result (\ref{eq4.8}) is equivalent to (\ref
{eq2.2}), and so (\ref{eq4.2}),
representing (\ref{eq2.2}) in the case $k=2$, implies (\ref{eq2.2})
for general $k\geq
2$, as had to be shown.

\textit{Step 2: Proof of (\ref{eq2.2}) when $k=2$.}
In this section, we shall show that, if $k=2$ and (\ref{eq2.5}) holds
for $\rho_{12}=e_2 n_1/(e_1 n_2)$, then (\ref{eq2.2}) holds.

To this end, define
\[
T=(S-ES)/(\var S)^{1/2} =\frac{{\bar X}_1+{\bar X}_2-\mu_1-\mu_2}{(n_1^{-1}\si_1^2+n_2^{-1}\si_2^2)^{1/2}} =c_1
T_1+c_2 T_2,
\]
where $T_j=({\bar X}_j-\mu_j)/(n_j^{-1}\si_j^2)^{1/2}$ and $c_1$ and $c_2$
are defined as at (\ref{eq2.9}). In this notation,
%
%
\begin{eqnarray}
\label{eq4.9} P(T\leq x)&=&P(c_1 T_1+c_2
T_2\leq x) =E \bigl\{P(c_1 T_1\leq
x-c_2 T_2\mid T_2) \bigr\}
\nonumber
\\
&=&E \biggl\{\Phi \biggl(\frac{x-c_2 T_2}{ c_1} \biggr) +n_1^{-1/2}A_1
\biggl(\frac{x-c_2 T_2}{ c_1} \biggr) +n_1^{-1/2}D_1
\biggl(\frac{x-c_2 T_2}{ c_1} \biggr) \biggr\}
\nonumber
\\
&&{} + %
\cases{ o \bigl(n^{-1/2} \bigr)& if $
\max_j E|X_{j1}|^3<\infty$
\cr
O
\bigl(n^{-1} \bigr)& if $\max_j E|X_{j1}|^4<
\infty$,} %
\end{eqnarray}
where $A_j$ and $D_j$ will refer to the smooth and discontinuous terms,
respectively, in the $n_j^{-1/2}$ component of an Edgeworth expansion of
the distribution of $T_j$ for $j=1,2$. In particular, $n_j^{-1/2}A_j$ and
$n_j^{-1/2}D_j$ are the counterparts of the second and third terms,
respectively, on the right-hand side of formula (35) p. 56 of Esseen
\cite{Ess45}.

Writing $B$ for either $\Phi$ or $A_1$, appearing on the right-hand
side of (\ref{eq4.9}), we have:
\begin{eqnarray*}
E \biggl\{B \biggl(\frac{x-c_2 T_2}{ c_1} \biggr) \biggr\} &=&\int B \biggl(
\frac{x-c_2 u}{ c_1} \biggr) \,dP(T_2\leq u)
\\
&=&\frac{c_2}{ c_1} \int B' \biggl(\frac{x-c_2 u}{ c_1} \biggr)
P(T_2\leq u) \,du.
\end{eqnarray*}
As in the argument leading to (\ref{eq4.6}) it can be shown that the
discontinuous term $n_2^{-1/2}D_2$, in the Edgeworth expansion of
$P(T_2\leq x)$, contributes only $O(n^{-1})$. Therefore, if we write
$\Ed_{2}(u)$ for the Edgeworth approximation to $P(T_2\leq u)$ that
includes the leading Gaussian term, plus the continuous part of the
component of order $n_2^{-1/2}$, and neglects everything else, we deduce
from (\ref{eq4.9}) that
%
%
\begin{eqnarray}
\label{eq4.10} P(T\leq x)&=&\int \biggl\{\Phi \biggl(\frac{x-c_2 u}{ c_1} \biggr)
+n_1^{-1/2}A_1 \biggl(\frac{x-c_2 u}{ c_1}
\biggr) \biggr\} \,d_u\Ed_{2}(u)
\nonumber
\\
&&{} +n_1^{-1/2}E \biggl\{D_1 \biggl(
\frac{x-c_2 T_2}{ c_1} \biggr) \biggr\} + %
\cases{ o \bigl(n^{-1/2}
\bigr)& if $\max_j E|X_{j1}|^3<\infty$
\cr
O
\bigl(n^{-1} \bigr)& if $\max_j E|X_{j1}|^4<
\infty$. }\qquad  %
\end{eqnarray}

Now we turn our attention to:
%
%
\begin{eqnarray}
\label{eq4.11}E \biggl\{D_1 \biggl(\frac{x-c_2 T_2}{ c_1} \biggr)
\biggr\} &=&\int D_1 \biggl(\frac{x-c_2 u}{ c_1} \biggr)
\,dP(T_2\leq u)
\nonumber
\\
&=&-\int P(T_2\leq u) \,d_uD_1 \biggl(
\frac{x-c_2 u}{ c_1} \biggr)
\\
&=&I_1(x)+n_2^{-1/2}I_2(x) +
\cases{ o(1) & if $\max_j E|X_{j1}|^3<
\infty$
\cr
O \bigl(n^{-1/2} \bigr)& if $\max_j
E|X_{j1}|^4<\infty$, } %
\nonumber
\\
\label{eq4.12}&=&n_2^{-1/2}I_2(x) + %
\cases{ o(1)& if $\max_j E|X_{j1}|^3<\infty$
\cr
O \bigl(n^{-1/2} \bigr)& if $\max_j E|X_{j1}|^4<
\infty$, } %
\end{eqnarray}
where
%
%
\begin{equation}
\label{eq4.13} I_1(x)=\int D_1 \biggl(
\frac{x-c_2 u}{ c_1} \biggr) \phi(u) \,du, \qquad I_2(x)=\int
D_1 \biggl(\frac{x-c_2 u}{c_1} \biggr) \,dD_2(u).
\end{equation}
To obtain the third identity in the string of formulae leading to
(\ref{eq4.12}), we used the integration by parts step at (\ref
{eq4.11}), a short Taylor
expansion of $P(T_2\leq u)$ with a remainder of $o(n^{-1/2})$ if $\max_j
E|X_{j1}|^3<\infty$ and $O (n^{-1})$ if $\max_j E|X_{j1}|^4$, and
the fact that $\int|dD_1|=O(n^{1/2})$ uniformly in $x$. (This can be
deduced either directly or by making use of (\ref{eq4.14}) below.)
Finally, it can be shown, arguing as in the proof of (\ref{eq4.6}), that
$I_1(x)=O(n^{-1/2})$, from which (\ref{eq4.12}) follows.

Note too that, with $\si_j$ defined as immediately above (\ref{eq2.1}),
%
%
\begin{eqnarray}
\label{eq4.14} D_j(x)&=&\frac{e_j}{\si_j} \psi \biggl\{
\frac{(x-\xi_{1n}) \si_j n_j^{1/2}}{ e_j} \biggr\} \phi (x)
\nonumber
\\
&=&\frac{e_j}{\si_j} \psi \biggl[\frac{\si_j n_j^{1/2}x}{ e_j} - \bigl\{(n_j
x_j/e_j)-\lfloor n_j x_j/e_j
\rfloor \bigr\} \biggr] \phi(x),
\end{eqnarray}
where, as in Sections~\ref{sec1} and~\ref{sec2}, $\psi(x)=\lfloor
x\rfloor-x+{\frac{1}{2}}$,
$\lfloor
x\rfloor$ is the largest integer not strictly exceeding $x$, and $\xi_{jn}=e_j (\si_j n_j^{1/2})^{-1}\{(n_j x_j/e_j)-\lfloor n_j
x_j/e_j\rfloor\}
$ if the lattice is located at points $x_j+\nu e_j$ for integers $\nu
$, see Esseen \cite{Ess45}, (29), (31) and (35) pp. 55/56. Defining
$\ga=(e_1e_2/\si_1\si_2)$, as in Section~\ref{sec2.3}; putting
\[
\psi_j(x)=\psi \biggl\{\frac{\si_j n_j^{1/2}}{ e_j} (x-\xi_{jn})
\biggr\} ;
\]
and noting that, by (\ref{eq4.14}), $D_j(x)=(e_j/\si_j) \psi_j(x)
\phi(x)$;
we deduce that
\begingroup
\abovedisplayskip=6.5pt
\belowdisplayskip=6.5pt
\begin{eqnarray*}
I_2(x)/\ga&=&\frac{1}{\ga}\int D_1 \biggl(
\frac{x-c_2 u}{ c_1} \biggr) \,dD_2(u) =\int(\psi_1 \phi)
\biggl(\frac{x-c_2 u}{ c_1} \biggr) \,d \bigl\{\psi_2(u) \phi(u) \bigr\}
\\
&=&\int(\psi_1 \phi) \biggl(\frac{x-c_2 u}{ c_1} \biggr) \bigl\{
\phi(u) \,d\psi_2(u)+\psi_2(u) \,d\phi(u) \bigr\}
\\
&=&\int(\psi_1 \phi) \biggl(\frac{x-c_2 u}{ c_1} \biggr) \phi(u) \, d
\psi_2(u) +{\frac{1}{2}}\int(\psi_1 \phi) \biggl(
\frac{x-c_2 u}{ c_1} \biggr) \psi_2(u) \,d\phi(u).
\end{eqnarray*}
The last-written integral equals $O(1)$, uniformly in $x$, and so, with
$I_2$ as at (\ref{eq4.13}),
%
%
\begin{equation}
\label{eq4.15} I_2(x)=\ga I_3(x)+O(1),
\end{equation}
uniformly in $x$, where
\[
I_3(x)=\int(\psi_1 \phi) \biggl(\frac{x-c_2 u}{ c_1}
\biggr) \phi(u) \,d\psi_2(u).
\]

Since $\psi_2$ has jumps of size $+1$ at points $u$ where $(u-\xi_{2n}) \si_2 n_2^{1/2}/e_2$ is an integer, i.e. $u=u_\nu\equiv\xi_{2n}+e_2 (\si_2 n_2^{1/2})^{-1}\nu$ for an integer $\nu$, then
%
%
\begin{eqnarray}
\label{eq4.16} I_3(x)&=&\sum_\nu(
\psi_1 \phi) \biggl(\frac{x-c_2 u_\nu}{ c_1} \biggr) \phi(u_\nu)
\nonumber
\\
&=&\sum_\nu\phi \biggl(\frac{x}{ c_1} -
\frac{c_2}{ c_1} \bigl\{\xi_{2n}+e_2 \bigl(
\si_2 n_2^{1/2} \bigr)^{-1} \nu \bigr\}
\biggr) \phi \bigl\{\xi_{2n}+e_2 \bigl(\si_2
n_2^{1/2} \bigr)^{-1}\nu \bigr\}
\nonumber
\\
&&{} \times \psi \biggl\{\xi_n(x)-\frac{e_2 n_1}{ e_1 n_2} \nu \biggr\},
\end{eqnarray}
where $\xi_n$ is as at (\ref{eq2.8}) and we have used the fact that
\begin{eqnarray*}
\psi_1 \biggl(\frac{x-c_2 u_\nu}{ c_1} \biggr) &=&\psi \biggl[
\frac{\si_1 n_1^{1/2}}{ e_1} \bigl\{(x-c_2 u_\nu)
c_1^{-1} -\xi_{1n} \bigr\} \biggr]
\\
&=&\psi \biggl\{\frac{(x-c_2 u_\nu-c_1 \xi_{1n}) \si_1 n_1^{1/2}}{ c_1
e_1} \biggr\}
\\
&=&\psi \biggl[\frac{\{x-(c_1 \xi_{1n}+c_2 \xi_{2n})-c_2 e_2 (\si_2
n_2^{1/2})^{-1}\nu\}
\si_1 n_1^{1/2}}{ c_1 e_1} \biggr]
\\
&=&\psi \biggl\{\xi_n(x)-\frac{e_2 n_1}{ e_1 n_2} \nu \biggr\},
\end{eqnarray*}
with
\[
\xi_n(x)= \bigl\{x-(c_1 \xi_{1n}+c_2
\xi_{2n}) \bigr\} \frac{\si_1 n_1^{1/2}}{
c_1 e_1}.\vadjust{\eject}
\]
\endgroup

Recall that $\xi_{jn}=e_j (\si_j n_j^{1/2})^{-1}\{(n_j
x_j/e_j)-\lfloor
n_j x_j/e_j\rfloor\}$ if the lattice is located at points $x_j+\nu e_j$
for integers $\nu$. In particular, $\xi_{jn}=O(n^{-1/2})$ for $j=1,2$.
Therefore, Taylor expanding the arguments of the functions $\phi$ at
(\ref{eq4.16}), and defining
%
%
\begin{equation}
\label{eq4.17} I_4(x)=\sum_\nu\phi
\biggl(\frac{x}{ c_1} -\frac{e_2 n_1^{1/2}}{\si_1 n_2} \nu \biggr) \phi \bigl
\{e_2 \bigl(\si_2 n_2^{1/2}
\bigr)^{-1}\nu \bigr\} \psi \biggl\{\xi_n(x)-
\frac{e_2 n_1}{ e_1 n_2} \nu \biggr\},
\end{equation}
we deduce that
%
%
\begin{equation}
\label{eq4.18} I_3(x)=I_4(x)+O(1),
\end{equation}
uniformly in $x$. Combining (\ref{eq4.10}), (\ref{eq4.12}), (\ref
{eq4.15}) and (\ref{eq4.18}), we
deduce that
%
%
\begin{eqnarray}
\label{eq4.19} P(T\leq x)&=&\int \biggl\{\Phi \biggl(\frac{x-c_2 u}{ c_1} \biggr)
+n_1^{-1/2}A_1 \biggl(\frac{x-c_2 u}{ c_1}
\biggr) \biggr\} \,d_u\Ed_{2}(u)
\nonumber
\\
&&{}+(n_1 n_2)^{-1/2}\ga I_4(x) +
\cases{ o \bigl(n^{-1/2} \bigr)& if $\max_j
E|X_{j1}|^3<\infty$
\cr
O \bigl(n^{-1} \bigr)& if
$\max_j E|X_{j1}|^4<\infty$.} %
\end{eqnarray}

If we can show that
%
%
\begin{equation}
\label{eq4.20} \sup_{-\infty<x<\infty} \bigl|I_4(x)\bigr|=o \bigl(n^{1/2}
\bigr)
\end{equation}
then it will follow from (\ref{eq4.19}), in cases where $\max_j
E|X_{j1}|^3<\infty$, that
%
%
\begin{equation}
\label{eq4.21} P(T\leq x)=\int \biggl\{\Phi \biggl(\frac{x-c_2 u}{ c_1} \biggr)
+n_1^{-1/2}A_1 \biggl(\frac{x-c_2 u}{ c_1}
\biggr) \biggr\} \,d_u\Ed_{2}(u) +o \bigl(n^{-1/2}
\bigr).
\end{equation}
The right-hand side here is Edgeworth expansion we would expect the
distribution of $T$ to enjoy if we were able to ignore the latticeness
of the distributions of $X_{j1}$ for $j=1,2$. That is, (\ref{eq4.21})
is just
(\ref{eq2.2}) in the particular case $k=2$. Therefore, provided (\ref
{eq4.20}) holds
then we shall have shown that (\ref{eq2.2}) holds whenever $k=2$. It
remains to
derive (\ref{eq4.20}).

\textit{Step 3: Proof of (\ref{eq4.20}).}
Given $\ep>0$, partition the set of all integers into adjacent blocks
$\cN_\ell$, for $-\infty<\ell<\infty$, where each block consists
of just $2 \lfloor n^{1/2}\ep\rfloor+1$ consecutive integers, and the
central integer is denoted by ${\bar\nu}_\ell$. Recalling the definition
of $I_4(x)$ at (\ref{eq4.17}), we deduce that
%
%
\begin{equation}
\label{eq4.22} I_4=\sum_{-\infty<\ell<\infty}
J_{1,\ell},
\end{equation}
where
%
%
\begin{equation}
\label{eq4.23} J_{1,\ell}(x)=\sum_{\nu\in\cN_\ell} \phi
\biggl(\frac{x}{ c_1} -\frac{e_2 n_1^{1/2}}{\si_1 n_2} \nu \biggr) \phi \bigl
\{e_2 \bigl(\si_2 n_2^{1/2}
\bigr)^{-1}\nu \bigr\} \psi \biggl\{\xi_n(x)-
\frac{e_2 n_1}{ e_1 n_2} \nu \biggr\}.
\end{equation}
Now,
%
%
\begin{equation}
\label{eq4.24} J_{1,\ell}=J_{2,\ell}+R_\ell,
\end{equation}
where
%
%
\begingroup
\abovedisplayskip=6.5pt
\belowdisplayskip=6.5pt
\begin{equation}
\label{eq4.25} J_{2,\ell}(x)=\phi \biggl(\frac{x}{ c_1} -
\frac{e_2 n_1^{1/2}}{\si_1 n_2} {\bar\nu}_\ell \biggr) \phi \bigl\{e_2
\bigl(\si_2 n_2^{1/2} \bigr)^{-1}{\bar
\nu}_\ell \bigr\} \sum_{\nu\in\cN_\ell} \psi \biggl\{
\xi_n(x)-\frac{e_2 n_1}{ e_1 n_2} \nu \biggr\}
\end{equation}
and $R_\ell$ is defined naively by (\ref{eq4.24}). Given an integer
$r$, let
$\ell(r)$ denote the unique value of $\ell$ such that $r\in\cN_\ell
$. Then, since $|\psi|\leq1$,
%
%
\begin{eqnarray}
\label{eq4.26} \biggl|\sum_{-\infty<\ell<\infty} R_\ell \biggl| &
\leq& \sum_r \biggr|\phi \biggl(\frac{x}{ c_1} -
\frac{e_2 n_1^{1/2}}{\si_1 n_2} r \biggr) \phi \bigl\{e_2 \bigl(\si_2
n_2^{1/2} \bigr)^{-1}r \bigr\}
\nonumber
\\
&&\hspace*{50pt} {} -\phi \biggl(\frac{x}{ c_1} -\frac{e_2 n_1^{1/2}}{\si_1 n_2} {\bar
\nu}_{\ell(r)} \biggr) \phi \bigl\{e_2 \bigl(\si_2
n_2^{1/2} \bigr)^{-1}{\bar\nu}_{\ell
(r)}
\bigr\} \biggr|
\nonumber
\\
&\leq& C_1 \ep n^{1/2},
\end{eqnarray}
where the constant $C_1$ does not depend on $\ep$ or $n$.

Let $\rho=e_2 n_1/(e_1 n_2)$, and define
\[
\chi_\cN(z,\rho)\equiv\sum_{\nu\in\cN}
\psi(z-\rho\nu).
\]
In this notation,
%
%
\begin{equation}
\label{eq4.27} J_{2,\ell}=\phi \biggl(\frac{x}{ c_1} -
\frac{e_2 n_1^{1/2}}{\si_1 n_2} {\bar\nu}_\ell \biggr) \phi \bigl\{e_2
\bigl(\si_2 n_2^{1/2} \bigr)^{-1}{\bar
\nu}_\ell \bigr\} \chi_{\cN_\ell} \bigl\{\xi_n(x),
\rho \bigr\}.
\end{equation}

If we can prove that, whenever the set $\cN$ consists of $|\cN|$
consecutive integers and $C_2<C_3$ are positive constants,
%
%
\begin{equation}
\label{eq4.28} \sup_{C_2 n^{1/2}\leq|\cN|\leq C_3 n^{1/2}} \sup_{-\infty<z<\infty} \bigl|\chi_\cN(z,
\rho)\bigr|=o \bigl(n^{1/2} \bigr)
\end{equation}
as $|\cN|\rai$, then it will follow from (\ref{eq4.27}) that
%
%
\begin{eqnarray}
\label{eq4.29} \biggl|\sum_{-\infty<\ell<\infty} J_{2,\ell} \biggr| &=&o
\biggl[n^{1/2}\sum_{-\infty<\ell<\infty} \phi \biggl(
\frac{x}{ c_1} -\frac{e_2 n_1^{1/2}}{\si_1 n_2} {\bar\nu}_\ell \biggr) \phi
\bigl\{e_2 \bigl(\si_2 n_2^{1/2}
\bigr)^{-1}{\bar\nu}_\ell \bigr\} \biggr]
\nonumber
\\
&=&o \bigl(n^{1/2} \bigr),
\end{eqnarray}
\endgroup
for each $\ep>0$, since the series on the first right-hand side of
(\ref{eq4.29}) is bounded uniformly in $n$. (To appreciate why,
observe that
${\bar\nu}_\ell$ is approximately an integer multiple of $n^{1/2}$,
plus a
constant.) Note that, since the left-hand side of (\ref{eq4.28})
involves the
supremum over~$z$, then that quantity does not depend on the location
of the set $\cN$ on the line, only on the number of consecutive
integers it contains.

The desired result (\ref{eq4.20}) follows from (\ref{eq4.22}), (\ref
{eq4.24}), the fact that
(\ref{eq4.26}) holds for each $\ep>0$, and (\ref{eq4.29}). To
complete the proof of
(\ref{eq4.20}), we shall derive (\ref{eq4.28}). Specifically, we
shall prove that, in
cases where (\ref{eq2.5}) is satisfied for $\rho_{12}=\rho=e_2 n_1/(e_1
n_2)$, (\ref{eq4.28}) obtains.

Assume that $\cN$ consists of $p$ consecutive integers, where $C_2
n^{1/2}\leq p\leq C_3 n^{1/2}$. Koksma's inequality (see, e.g.,
Theorems 1.3 and 5.1, pp. 91 and 143 of \cite{Kui74}), and the Erd\H
os-Tur\'an inequality (see, e.g., formula (2.42), p. 114 of \cite
{Kui74}), can be combined to prove that, for all integers $m\geq1$,
%
%
\begin{eqnarray}
\label{eq4.30} \sup_z \bigl|\chi_\cN(z,\rho_{12})\bigr|
&\leq& C_4 \Biggl\{\frac{p}{ m} +\sum
_{\ell=1}^m \frac{1}{\ell} \sup_z \Biggl|
\sum_{r=1}^p \exp(2\pi i \ell r
\rho_{12})\Biggr | \Biggr\}
\nonumber
\\
&\leq& C_4 \Biggl\{\frac{p}{ m} +\sum
_{\ell=1}^m \frac{1}{\ell|\sin(\ell\rho_{12}\pi)|} \Biggr\} ,
\end{eqnarray}
where $C_4$ is an absolute constant. Since (\ref{eq2.5}) is assumed to hold
with $(j_1,j_2)=(1,2)$ then, for each fixed $m$,
\[
\max_{1\leq\ell\leq m} \bigl|\sin(\ell\rho_{12}\pi)\bigr|^{-1}=o
\bigl(n^{1/2} \bigr).
\]
Hence, by (\ref{eq4.30}),
%
%
\begin{equation}
\label{eq4.31} \sup_z \bigl|\chi_\cN(z,\rho_{12})\bigr|
\leq\frac{C_3 C_4 n^{1/2}}{ m} +o \bigl(n^{1/2} \bigr),
\end{equation}
where the $o(n^{1/2})$ term is of that order uniformly in $\cN$ such
that $C_2 n^{1/2}\leq|\cN|\leq C_3 n^{1/2}$. However, $m$ can be
taken arbitrarily large, and none of $C_2$, $C_3$ and $C_4$ depends on
$m$ or~$n$. Therefore, (\ref{eq4.31}) implies (\ref{eq4.28}).

\subsubsection{\texorpdfstring{Proof of part (ii) of Theorem~\protect\ref{the1}}
{Proof of part (ii) of Theorem 1}}\label{sec4.1.2}

We can write
\[
{\bar X}_1+\cdots+{\bar X}_k=\frac{e_1}{ n_1}
(Y_1+\cdots+Y_k)+\mu,
\]
where $\mu$ is deterministic and, for each $j$, $Y_j$ is the sum of
$n_j$ random variables $Y_{j1},\ldots,Y_{jn_j}$, each having a lattice
distribution (not depending on $n$) supported on the set of points
$\rho_{1j} \ell$ for $\ell\in\ZZ$, and with the $Y_{ji}$s being
totally independent. Of course, $\rho_{11}=1$. Since each $\rho_{j_1j_2}$ equals a rational number, not depending on $n$, then the set
$\bigcup_j \{\rho_{1j} \ell, \ell\in\ZZ\}$ can itself be represented
as a maximal lattice, $\cL$ say, not depending on $n$. The
distribution of
\[
Y_1+\cdots+Y_k=\sum_{j=1}^k
\sum_{i=1}^{n_j} Y_{ji}
\]
can be viewed as the distribution of the sum of $n=n_1+\cdots+n_k$
independent and identically distributed random variables each having a
mixture distribution, $D_n$ say, with support confined to~$\cL$.
Although $D_n$ depends on $n$, since it is always supported on the same
lattice, standard methods can be used to derive an Edgeworth expansion
of the distribution of $Y_1+\cdots+Y_k$, from which it can be seen
that there is a nonvanishing discontinuous term, not present in (\ref{eq2.2}).

\subsection{\texorpdfstring{Proof of Theorem~\protect\ref{the2}}{Proof of Theorem 2}}\label{sec4.2}

\textit{Step 1: Proof that it is sufficient to consider the case $k=2$.}
We give the argument only in outline, since it parallels that in step 1
of the derivation of Theorem~\ref{the1}. Suppose it is possible to
derive the
version of (\ref{eq4.1}) where the remainder $o(n^{-1/2})$ is replaced by
$O(n^{\xi-1})$, for all $\xi>0$. Then, as in the earlier proof, we
have (\ref{eq4.2}) where the remainder term is $O(n^{\xi-1})$, for
all $\xi
>0$, instead of $o(n^{-1/2})$. The string of arguments leading to (\ref{eq4.5})
holds without change, as too does (\ref{eq4.6}). Combining the revised
(\ref{eq4.2})
with the old (\ref{eq4.5}) and (\ref{eq4.6}) we deduce the following
version of (\ref{eq4.8}):
\begin{eqnarray*}
P(S\leq x)&=&\int \biggl(\Phi \biggl\{\frac{x-r}{(\var
S_1)^{1/2}} \biggr\} +
\frac{\be_1}{6 n^{1/2}} \biggl[1- \biggl\{\frac{x-r}{(\var S_1)^{1/2}} \biggr\}^{ 2}
\biggr] \phi \biggl\{\frac{x-r}{(\var S_1)^{1/2}} \biggr\} \biggr)
\\
&&{} \times d \biggl\{\Phi(x/\tau_2)+n^{-1/2}{
\frac{1}{6}}\be_2 \bigl\{1-(r/\tau_2)^2
\bigr\} \phi(r/\tau_2) \biggr\} +O \bigl(n^{\xi-1} \bigr),
\end{eqnarray*}
uniformly in $x$ and for all $\xi>0$. This formula is equivalent to
(\ref{eq2.2}), with the remainder there replaced by $O(n^{\xi-1})$,
and so we
have shown that it suffices to consider $k=2$.

\textit{Step 2: Completion of proof of Theorem~\ref{the2}.}
Combining (\ref{eq4.10}) and (\ref{eq4.12}) in the case $\max_j
E|X_{j1}|^4<\infty$,
and noting (\ref{eq4.15}) and (\ref{eq4.18}), we deduce the version
of (\ref{eq4.19}) when
$\max_j E|X_{j1}|^4<\infty$.

Next we reintroduce the notation noted below (\ref{eq2.8}), where $\al
\in
(0,{\frac{1}{2}})$, $\cN_\ell$ (for $-\infty<\ell<\infty$) is a partition
of the set of all integers into adjacent blocks each containing $2
\lfloor
n^\al\rfloor+1$ consecutive integers, ${\bar\nu}_\ell$ is the central
integer in $\cN_\ell$, and $\nu_\ell=\nu-{\bar\nu}_\ell$ for
$\nu\in
\cN_\ell$. Property (\ref{eq4.22}) continues to hold, with
$J_{1,\ell}$
still given by (\ref{eq4.23}). Again we define $R_\ell$ and
$J_{2,\ell}$ by
(\ref{eq4.24}) and (\ref{eq4.25}). However, this time we give an
expansion for, rather
than an upper bound to, $R_\ell$. As a first step, note that
\begin{eqnarray*}
R_\ell(x)&=&J_{1,\ell}(x)-J_{2,\ell}(x)
\\
&=&\sum_{\nu\in\cN_\ell} \biggl[\phi \biggl\{\frac{x}{ c_1}
-\frac{e_2 n_1^{1/2}}{\si_1 n_2} ({\bar\nu}_\ell+\nu_\ell) \biggr\} \phi
\bigl\{e_2 \bigl(\si_2 n_2^{1/2}
\bigr)^{-1}({\bar\nu}_\ell+\nu_\ell ) \bigr\}
\\
&&{} -\phi \biggl(\frac{x}{ c_1} -\frac{e_2 n_1^{1/2}}{\si_1 n_2} {\bar
\nu}_\ell \biggr) \phi \bigl\{e_2 \bigl(\si_2
n_2^{1/2} \bigr)^{-1}{\bar\nu}_\ell
\bigr\} \biggr] \psi \biggl\{\xi_n(x)-\frac{e_2 n_1}{ e_1 n_2} \nu \biggr
\}.
\end{eqnarray*}
Taylor-expanding, and using the argument in the paragraph immediately
below that containing (\ref{eq2.10}), we deduce that
%
%
\begin{eqnarray}
\label{eq4.32} &&\sum_{-\infty<\ell<\infty} R_\ell(x) =
\sum_{r=1}^{r_0} \sum
_{-\infty<\ell<\infty} \frac{\phi_r ({\bar\nu}_\ell /n_1^{1/2},x )}{ r! n_1^{r/2}} \sum_{\nu\in\cN_\ell}
\nu_\ell^r \psi \biggl\{\xi_n(x)-
\frac{e_2 n_1}{ e_1 n_2} \nu \biggr\}
\nonumber
\\
&&\hspace*{170pt} {}+O \bigl(n^\al\cdot n^{1/2}\cdot
n^{(r_0+1) \{\al
-(1/2)\}} \bigr),\quad
\end{eqnarray}
uniformly in $x$. Adding $\sum_\ell J_{2,\ell}$ to either side of
(\ref{eq4.32}) has the effect, on the right-hand side, of changing the
range of
summation of the first series to $0\leq r\leq r_0$. Therefore,
%
%
\begin{equation}
\label{eq4.33} \sum_{-\infty<\ell<\infty} \bigl\{J_{2,\ell}(x)+R_\ell(x)
\bigr\} =\ga^{-1}K_n(x) +O \bigl(n^\al\cdot
n^{1/2}\cdot n^{(r_0+1) \{\al-(1/2)\}} \bigr) ,
\end{equation}
uniformly in $x$, where $\ga=\prod_{j=1,2} (e_j/\si_j)$ and $K_n$ is
at (\ref{eq2.10}). If $r_0\geq4\al/(1-2\al)$, as stipulated in
Theorem~\ref{the2},
then the ``$O$'' remainder in (\ref{eq4.33}) is just $O(1)$. In this case,
%
%
\begin{equation}
\label{eq4.34} \ga I_4(x)=\ga\sum_{-\infty<\ell<\infty}
J_{1,\ell}(x) =\ga\sum_{-\infty<\ell<\infty} \bigl
\{J_{2,\ell}(x)+R_\ell(x) \bigr\} =K_n(x)+O(1),
\end{equation}
uniformly in $x$. Part (i) of Theorem~\ref{the2}, which addresses only
the case
$k=2$, follows from (\ref{eq4.19}) and (\ref{eq4.34}). Part (ii) of
Theorem~\ref{the2}, in the
case $k=2$, follows from (\ref{eq4.34}) and (\ref{eq2.36}). In view
of Part 1 of the
proof of Theorem~\ref{the2}, this is sufficient to complete the proof
of the theorem.

\section*{Acknowledgements}
We are grateful to Professor Roger Heath-Brown for helpful discussion.
The research was supported by the Australian Research Council and the
National Science Foundation.


%

%

\end{document}